\documentclass[12pt]{amsart}
\usepackage{amssymb,latexsym}
\usepackage{enumerate}
\usepackage{amscd}
\usepackage{verbatim}

\theoremstyle{plain}
\newtheorem{theorem}{Theorem}
\newtheorem{corollary}[theorem]{Corollary}
\newtheorem{lemma}[theorem]{Lemma}

\theoremstyle{definition}

\theoremstyle{remark}
\newtheorem*{remark}{Remark}
\newtheorem*{notation}{Notation}
\newtheorem{example}{Example}
\newtheorem*{catnote}{Note for category theorists}

\numberwithin{equation}{section}
\numberwithin{theorem}{section}
\numberwithin{example}{section}

\def\pair#1#2{\langle #1, #2\rangle}

\def\semitimes{\ltimes}

\def\Mof#1{\Cal M(A)}

\newcommand{\N}{{\mathbb N}}
\newcommand{\Z}{{\mathbb Z}}

\def\Set{{\text{\bf Set}}}

\def\Grp{{\text{\bf Grp}}}

\def\CRng{{\text{\bf CRng}}}
\def\Ab{{\text{\bf Ab}}}
\def\Pnt{{\text{\bf Pnt}}}
\def\Alg{{\text{\bf Alg}}}

\def\Ov{{\text{\bf Ov}}}

\def\Mod{{\text{\bf Mod}}}

\def\natur{{\overset{\scriptscriptstyle\bullet}{\to}}}

\newcommand{\Con}{\operatorname{Con}}

\newcommand{\End}{\operatorname{End}}

\newcommand{\Pol}{\operatorname{Pol}}

\newcommand{\Sub}{\operatorname{Sub}}

\newcommand{\Ke}{\operatorname{Ke}}
\newcommand{\Co}{\operatorname{Co}}

\newcommand{\MyIm}{\operatorname{Im}}
\newcommand{\nat}{\operatorname{nat}}

\newcommand{\Clo}{\operatorname{Clo}}

\newcommand{\I}{\operatorname{I}}

\def\Var{{\mathbf V}}

\def\lsil{\lbrack\!\lbrack}
\def\rsil{\rbrack\!\rbrack}

\def\congruence{on\-gru\-ence\discretionary{-}{}{-}}
\def\conM/{c\congruence mod\-u\-lar}
\def\ConM/{C\congruence mod\-u\-lar}
\def\conD/{c\congruence dis\-trib\-u\-tive}
\def\ConD/{C\congruence dis\-trib\-u\-tive}
\def\conP/{c\congruence per\-mut\-a\-ble}
\def\ConP/{C\congruence per\-mut\-a\-ble}
\def\conMity/{\conM/\-i\-ty}
\def\ConMity/{\ConM/\-i\-ty}
\def\conDity/{c\congruence dis\-trib\-u\-tiv\-i\-ty}
\def\ConDity/{C\congruence dis\-trib\-u\-tiv\-i\-ty}
\def\conPity/{c\congruence per\-mut\-a\-bil\-i\-ty}
\def\ConPity/{C\congruence per\-mut\-a\-bil\-i\-ty}
\def\usprv/{un\-der\-ly\-ing-set-pre\-ser\-ving}

\def\ie/{{i.e.}}
\def\Ie/{{I.e.}}
\def\eg/{{e.g.}}
\def\Eg/{{E.g.}}
\def\etc/{{etc.}}

\newdimen\mysubdimen
\newbox\mysubbox

\def\subwhat#1#2#3{{
\setbox\mysubbox=\hbox{#3}
\mysubdimen=\wd\mysubbox
\setbox\mysubbox=\hbox{$#1#2$}
\ifnum\mysubdimen>\wd\mysubbox
\vtop{
\hbox to\mysubdimen{\hfil\box\mysubbox\hfil}
\nointerlineskip
\hbox{#3}}
\else
\mysubdimen=\wd\mysubbox
\vtop{
\box\mysubbox
\nointerlineskip
\hbox to\mysubdimen{#3}}
\fi
}}

\begin{document}

\title[Modulization and the Enveloping Ringoid]
{Modulization and the Enveloping Ringoid}
\author{William H. Rowan}
\address{PO Box 20161 \\
         Oakland, California 94620}
\email{whrowan@member.ams.org}

\date{\today}

\begin{abstract}
Let $A$ be an algebra in a variety $\mathbf V$. We study
the modulization of a pointed $A$-overalgebra $P$, show
that it is totally in any variety that $P$ is totally
in, and apply this theory to the
construction of the enveloping ringoid $\Z[A,\mathbf V]$.
\end{abstract}

\maketitle

\begin{table}
\begin{center}
\begin{tabular}{ | l | c | l | }
\hline
\quad $A\in\Var$ &
$_e\Z[A,\Var]_e$ & \quad Usual Name\\ \hline
group $G$ & $\Z[G]$ & group ring\\ \hline
abelian group $G$ & $\Z$ & the integers\\
\hline
group $G$ in a variety of groups $\Var_S$ &
$\Z_S[G]$ & universal $\Var_S$-envelope \cite{Kn66}\\
\hline
Lie algebra $L$ & $U(L)$ & universal enveloping algebra\\
\hline
Jordan algebra $J$ & $U(J)^{\text{op}}$
& $(\text{universal multiplication
envelope})^{\text{op}}$\\ \hline
commutative ring $R$ & $R$ & $R$ itself\\ \hline
ring $R$ & $R\otimes_{\Z} R^{\text{op}}$
& enveloping ring\\ \hline
\end{tabular}
\end{center}
\caption{Familiar Rings Equivalent to
Enveloping Ringoids}
\end{table}

\section*{Introduction}

Let $A$ be an algebra of some type $\Omega$. The concept
of an \emph{$A$-module} was first studied in the context
of category theory \cite{E}, and then generalized
\cite{B} in the form of what has come to be known as
\emph{Beck modules}. Our own work on this subject
\cite{R92,R96} includes the definition of $\Ab[A,\mathbf
V]$, the \emph{category of $A$-modules
totally in the variety $\mathbf V$}. This category is
equivalent to the earlier-defined category of Beck
modules, with some features which we feel are
advantageous.

Our previous
work also included the construction of the enveloping
ringoid $\Z[A,\mathbf V]$, over which the left modules
form a category isomorphic to $\Ab[A,\mathbf V]$.
Table~1 shows, in the
first column, an algebra in a familiar variety, and in the
second and third columns, a familiar ring equivalent to
the enveloping ringoid $\Z[A,\mathbf V]$. 
All of these
examples are of varieties of algebras with forgetful
functors to the variety of groups, although
there are additional examples
\cite{R92,R96} which are not of this type. When a variety
does have a forgetful functor to the variety of groups,
the enveloping ringoid $\Z[A,\mathbf V]$ is always
equivalent to a ring, and the equivalence is
coherent as $A$ ranges through $\mathbf V$
\cite{R96}.

In \cite{R96}, an approach to the definition of the
enveloping ringoid was given which used the concept
of \emph{modulization}, in which an object called a
\emph{pointed $A$-overalgebra} is made into an $A$-module.
The discussion in \cite{R96} was limited to algebras in
congruence-modular varieties
$\mathbf V$. The same construction was used more
generally, but without being explicit about applying the
modulization functor, in \cite{R92}.  In this paper, we
study modulization in detail for the general (not
necessarily congruence-modular) case, and use
modulization to construct the enveloping ringoid. In the
process, we give a better-motivated construction of the
enveloping ringoid than has been given previously.

The purposes of the paper also include describing a
calculus of polynomials, with its connections to the
theory of pointed $A$-overalgebras and $A$-modules, and
proving the simplifying fact that if a pointed overalgebra
is totally in a variety $\mathbf V$, then so is its
modulization.  Also, a detailed construction of the
enveloping ringoid has been given only in \cite{R92}, and
there is some value in publishing it, along with a proof
that the categories $\Z[A,\mathbf V]$-$\Mod$ and
$\Ab[A,\mathbf V]$ are not just equivalent, but
isomorphic. This has been a point of some confusion.

The plan of the paper is as follows. First there is a
section of preliminaries, including discussions of
clones, and of ringoids and their modules. Then, in
\S\ref{S:Modules} to \S\ref{S:SubsandQuotients}, we define
$A$-modules and pointed $A$-overalgebras, and
introduce related concepts.  Next, in
\S\ref{S:Polynomials}, we discuss the clone
$\Pol(A,\mathbf V)$ of
\emph{polynomials with coefficients in $A$}, and its
relationship to modules and pointed overalgebras.
Following that, in \S\ref{S:ZsubM}, we define the ringoid
$Z_M$ corresponding to an $A$-module $M$, which will later
be shown to be a canonical image of the enveloping
ringoid of $A$. In
\S\ref{S:Modulization}, we study the modulization
functor, and show that the modulization of a pointed
$A$-overalgebra totally in $\mathbf V$ is also totally in
$\mathbf V$. In \S\ref{S:Enveloping} and
\S\ref{S:Isomorphism}, we study
the enveloping ringoid, and show that the 
categories $\Z[A,\mathbf V]$-$\Mod$ and $\Ab[A,\mathbf V]$
are isomorphic. In \S\ref{S:Canonical}, we define the
canonical ringoid homomorphism from the enveloping ringoid
to $Z_M$.
In \S\ref{S:JandR} to \S\ref{S:Previous}, we prove that
this way of constructing the enveloping ringoid is
actually identical to one of the ways given in
\cite{R92}, although the discussion there is not as well
motivated. Only in
\S\ref{S:CM} do we consider the hypothesis of
congruence-modularity of the variety
$\mathbf V$, in sketching the simplifications in
the theory of the enveloping ringoid which become manifest
in that case.

\section*{Preliminaries}

\subsection{Category theory} We follow
\cite{MacL} in terminology and notation.

\subsection{Lattice theory} We will denote the least and
greatest elements of a lattice by $\bot$ and $\top$, and
if $L$ is a lattice, and $a$, $b\in L$, the interval
sublattice from $a$ to $b$ by $\I_L[a,b]$.

\subsection{Universal algebra} The basic
definitions of Universal Algebra can be found, for
example, in
\cite{B-S}. However, we have some preferences in
notation and in the way the subject is developed, as
follows:

In the definition of
an algebra, we prefer to allow an algebra to be empty.

If $A$ is an algebra, we denote its underlying set by
$|A|$.

Recall that a \emph{type} is an $\N$-tuple of sets, where
$\N$ stands for the set of natural numbers. If
$\Omega$ is a type, then we call the elements of
$\Omega_n$ the
\emph{$n$-ary elements} of $\Omega$.

If $\mathbf V$ is a variety of algebras of some type, we
will typically make no distinction between the category of
algebras in $\mathbf V$, and $\mathbf V$ itself.

$\mathcal F_\Omega S$ will stand for the \emph{free
algebra} of type $\Omega$ on the set $S$ (i.e., the
\emph{word algebra}), while
$\mathcal F_{\mathbf V}S$ will stand for its quotient, the
\emph{free algebra relative to the variety
$\mathbf V$}.

We define an $n$-ary term (of type $\Omega$) to be an
element of the algebra $T_n=\mathcal
F_\Omega\{\,x_1,\ldots,x_n\,\}$. If
$A$ is an algebra of type
$\Omega$, and $t$ is an
$n$-ary term of type $\Omega$, then $t^A$ will stand for
the $n$-ary function on $A$ that sends each $n$-tuple
$\mathbf a$ to the image of $t$, under the unique
homomorphism from $T_n$ to $A$
sending each $x_i$ to $a_i$.

\subsection{Clones}
\emph{Clones} are a part of Universal Algebra which will
be important to our development. A \emph{clone} $C$
consists of an
$\N$-tuple of sets
$C_n$, called the set of \emph{$n$-ary elements} of $C$,
together with some additional structure, as follows.
For each $n$, there are elements $\pi^C_{i,n}\in C_n$
as $i$ ranges from $1$ to $n$, called the
\emph{$i^{\text{th}}$ of $n$ projection}.  For each
$n$ and
$n'$, each $n'$-tuple of $n$-ary elements $\mathbf c$,
and each
$n'$-ary element $c'$, there is an $n$-ary element
$c'\mathbf c$, called the \emph{composite} of $c'$
and $\mathbf c$. The constants $\pi^C_{i,n}$ and clone
composition satisfy
\begin{enumerate}
\item[(C1)] For all $n$ and all $c\in C_n$,
$c\langle\pi^C_{1,n},\ldots,\pi^C_{n,n}\rangle=c$;
\item[(C2)] For all $n$ and $n'$, and all $i$ such that
$1\leq i\leq n$, and for each $n$-tuple $\mathbf c$ of
elements of
$C_{n'}$,
$\pi^C_{i,n}\mathbf c=c_i$; and
\item[(C3)] For all $n$, $n'$, and $n''$, if $\mathbf
c\in(C_n)^{n'}$, $\mathbf{c'}\in(C_{n'})^{n''}$, and
$c''\in C_{n''}$, we have
$c''(\mathbf{c'}\mathbf c)=(c''\mathbf c')\mathbf c$,
where $\mathbf{c'}\mathbf c$ stands for $\langle
c'_1\mathbf c,\ldots,c'_{n''}\mathbf c\rangle$.
\end{enumerate}

If $C$ and $C'$ are clones, then an $\N$-tuple of
mappings $\Phi_n:C_n\to C'_n$ is called a \emph{clone
homomorphism} if for all $n$ and $i$, we have
$\Phi_n(\pi^C_{i,n})=\pi^{C'}_{i,n}$ and for all $n$,
$n'$,
$\mathbf c\in C_n^{n'}$, and $c'\in C_{n'}$, we have
$\Phi_n(c'\mathbf c)=(\Phi_{n'}(c'))(\Phi_n(\mathbf c))$.

There are many natural examples of clones and clone
homomorphisms in Universal Algebra.
Given a set $S$, $\Clo S$ will denote the
\emph{clone of operations on $S$}. $\Clo_n S$ is the set
of $n$-ary operations on $S$, $\pi^{\Clo S}_{i,n}(\mathbf
s)=s_i$, and
given $n$, $n'$, an $n'$-tuple of $n$-ary operations
$\mathbf w$,
an
$n'$-ary operation
$w'$, and an
$n$-tuple $\mathbf s$ of elements of $S$, we define
$(w'\mathbf w)(\mathbf
s)=w'(w_1(\mathbf s),\ldots,w_{n'}(\mathbf
s))$. It is straightforward to verify that these
definitions make $\Clo S$ a clone.

If $\Omega$ is a type, we define
$\Clo\Omega$, the \emph{clone of terms of type
$\Omega$}, as follows: $\Clo_n\Omega=|T_n|$, the set
of
$n$-ary terms of type $\Omega$, defined above. We
define
$\pi^{\Clo\Omega}_{i,n}=x_i$ and, given $n$, $n'$,
an $n'$-tuple $\mathbf t$ of elements of $\Clo_n\Omega$,
and
$t'\in\Clo_{n'}\Omega$, we define $t'\mathbf t$ to
be $(t')^{T_n}(\mathbf t)$.
Note that $\Clo\Omega$ is a free clone on the type
$\Omega$; that is, there is a type homomorphism
$\Phi_\Omega:\Omega\to\Clo\Omega$ ($\N$-tuple of functions
$(\Phi_\Omega)_n:\Omega_n\to\Clo_n\Omega$) such that
$\pair{\Clo\Omega}{\Phi_\Omega}$ is a universal arrow from
the type $\Omega$ to the forgetful functor from clones to
types. $(\Phi_\Omega)_n$ sends $\omega\in\Omega_n$ to
$\omega(x_1,\ldots,x_n)$, and we have
$\omega^A=(\omega(\mathbf x))^A$ for every algebra $A$.

If $A$ is an algebra of type $\Omega$, the mappings
$\Phi^A_n:\Clo_n\Omega\to\Clo_n|A|$ given by $t\mapsto
t^A$ form a clone homomorphism. In fact,
specifying an algebra $A$ is the same as specifying
$S=|A|$, the underlying set, and a clone homomorphism
from $\Clo\Omega$ to $\Clo S$.

It is important to distinguish between terms and
\emph{term operations}, which are operations on the
underlying set of an algebra
$A$ of the form $t^A$ for some term $t$. (The term
operations on $A$ form a \emph{subclone} of $\Clo|A|$.)

If $\mathbf V$ is a variety of algebras of type
$\Omega$, then we may contrast $\Clo\Omega$ with
$\Clo\mathbf V$, the \emph{clone of the variety
$\mathbf V$}.  Elements of $\Clo_n\mathbf V$ can be
defined as equivalence classes of elements of
$\Clo_n\Omega$, under the equivalence relation that
relates terms
$t$, $t'$ such that $t^A=(t')^A$ for all $A\in\mathbf
V$. Or, we can define
$\Clo_n\mathbf V=\mathcal
F_{\mathbf V}\{\,x_1,\ldots,x_n\,\}$. In any case, there
is an associated clone homomorphism $\Phi_{\mathbf
V}:\Clo\Omega\to\Clo\mathbf V$. If $A$ is an algebra of
type $\Omega$, then we have $A\in\mathbf V$ iff
$\Phi^A$ factors through $\Phi_{\mathbf V}$.

For readability, when writing down unary terms, we
will write $x$ rather than
$x_1$.

\begin{catnote}
Clones are closely related to the
category-theoretic concept of \emph{theories}
\cite{L}. Given a clone $C$, we define a category $\bar
C$ with objects the natural numbers, and such that an
arrow from
$m$ to $n$ is an $n$-tuple of elements of $C_m$. We
define $1_n$ to be
$\langle\pi^C_{1,n},\ldots,\pi^C_{n,n}\rangle$. If
$\mathbf c\in\bar C(m,n)$ and $\mathbf{c'}\in\bar
C(n,p)$, then we define the composition of $\mathbf
c$ and $\mathbf{c'}$ to be $\mathbf{c'}\mathbf c$. Note
that conditions (C1), (C2), and (C3) translate very
simply to the
axioms for a category. $\bar C$ is not only a category,
but a theory, because each object $n$ is the $n$-fold
direct power of
$1$. If we construct the theory $\mathbf
T_{\mathbf V}$ corresponding to
$\Clo\mathbf V$ for some variety $\mathbf V$, then an
algebra in $\mathbf V$ can be thought of as a
product-preserving functor from $\mathbf T_{\mathbf V}$
to the category of sets.
\end{catnote}

\subsection{Commutator theory}
In \S\ref{S:CM}, which discusses special results applying
to the congruence-modular case, we will assume some
acquaintance with commutator theory for
congruence-modular varieties, as described in
\cite{F-McK}. In particular, we will use the notion of a
difference term.

\subsection{Ringoids, and their modules
and bimodules}

A \emph{ringoid} is a small additive category.  That is,
a ringoid is a small category $X$ such that the
hom-set
$X(a,b)$ between two objects $a$ and $b$ is an
abelian group, and composition is bilinear.
A \emph{ringoid homomorphism} is simply an additive
functor, i.e., a functor which is an abelian group
homomorphism on each hom-set.

If $X$ is a ringoid, a \emph{left $X$-module}
is an additive functor from $X$ to $\Ab$, the
category of abelian groups. Similarly, a \emph{right
$X$-module} is an additive functor from
$X^{\text{op}}$ to $\Ab$.  Homomorphisms of left or right
modules are simply natural transformations.  We write
$X$-$\Mod$ for the category of left $X$-modules, and
$\Mod$-$X$ for the category of right $X$-modules.

If $X$ is a ringoid, with set of objects $A$, we write
$_{a'}X_a$ rather than $X(a,a')$.  We write $1^X_a$, or
simply $1_a$, for the identity element of $_aX_a$, and
$_{a'}0^X_a$, or simply $_{a'}0_a$, for the zero element
of $_{a'}X_a$. We call the composition of arrows
\emph{multiplication}. For each $a\in A$,
$X_a$ will stand for the left $X$-module consisting of
the abelian groups
$_bX_a$, and $_aX$ will stand for the right $X$-module
consisting of the groups $_aX_b$. If
$Y$ is another ringoid and $f:X\to Y$ is a ringoid
homomorphism, we write $_{a'}f_a$ for the abelian group
homomorphism from
$_{a'}X_a$ to
$_{fa'}Y_{fa}$ defined by the additive functor $f$. If
$M$ is a left
$X$-module, we write $_aM$ rather than $Ma$, and for each
$m\in{_aM}$ and $x\in{_{a'}X_a}$, we write $xm$ rather
than $(Mx)(m)$. Finally, if
$M$ and $M'$ are left $X$-modules, and $\phi:M\natur M'$
is a homomorphism, then we simply write $\phi:M\to M'$,
and we write $_a\phi$ for the component of $\phi$ at $a$,
rather than $\phi a$.

Our notation allows us to express the category-theoretic
definitions of ringoid theory in a form more reminiscent
of ring theory. For example, suppose that $X$ is a
ringoid, with set of objects $A$, and that $M$ is a left
$X$-module. This means that $M$ is an $A$-tuple of
abelian groups, that for each $a$ and $m\in{_aM}$,
$1^X_am=m$, that $xm$ is $\Z$-bilinear in $x$ and in
$m$, and that for each
$a$,
$a'$,
$a''\in A$, each
$x\in{_aX_{a'}}$, and each $x'\in{_{a'}X_{a''}}$,
$(x'x)m=x'(xm)$.

As an example of a ringoid, let $\{\,_aM\,\}_{a\in A}$
be an
$A$-tuple of abelian groups for some index set $A$. We
define the ringoid
$\End(M)$ to have set of objects $A$, and for each $a$,
$a'\in A$, to have hom-set
$_{a'}\End(M)_a=\Ab({_aM},{_{a'}M})$. 
Note that $M$ is a left
$\End(M)$-module in an obvious way.

If $X$ is a ringoid, with set of objects $A$, and $J$
is an $A^2$-tuple of subgroups $_{a'}J_a$ of
the respective $_{a'}X_a$, then we say that $J$ is an {\it
ideal} of $X$ if for each $a$, $a'$, and $a''\in A$, and
each $x\in{_{a'}X_a}$ and $x'\in {_{a''}X_{a'}}$, if
either $x'\in{_{a''}J_{a'}}$ or
$x\in{_{a'}J_a}$, then $x'x\in{_{a''}J_a}$.

If $X$ and $X'$ are ringoids, $X$ has set of
objects $A$, and $f:X\to X'$ is a homomorphism, then we
define $\Ke f$, the \emph{kernel} of $f$, to be the
$A^2$-tuple of subgroups of the $_{a'}X_a$
defined by $_{a'}(\Ke f)_a=\Ke{_{a'}f_a}$. It is easy to
see that $\Ke f$ is an ideal of $X$.

If $X$ is a ringoid, having $A$ as set of
objects, and $J$ is an ideal of $X$, then we can define a
new ringoid $X/J$, again having $A$ as its set of
objects, by the formula
$_{a'}(X/J)_a={_{a'}X_a}/{_{a'}J_a}$, with identity
elements $1^{X/J}_a=1^X_a/{_aJ_a}$, and with the
multiplication well-defined by
$(x'/{_{a''}J_{a'}})(x/{_{a'}J_a})=x'x/{_{a''}J_a}$.  We
call $X/J$ the \emph{ quotient} of $X$ by the ideal $J$.
Clearly, the
$A^2$-tuple of mappings
${_{a'}X_a}\to{_{a'}(X/J)_a}$ defined by
$x\mapsto x/{_{a'}J_a}$, form a homomorphism with kernel
$J$, and we denote this homomorphism by $\nat J$.

To simplify notation, in what follows, we will typically
suppress the subscripts for $J$. Thus, we will write
$x/J$ rather than $x/{_{a'}J_a}$, and ${_{a'}X_a}/J$
rather than ${_{a'}X_a}/{_{a'}J_a}$.

If $X$ is a ringoid, with set of
objects $A$, then we call an $A^2$-tuple of subgroups
$_b{X'}_a\subseteq{_bX_a}$ a \emph{subringoid} if
$x\in{_bX_a}$, $x'\in{_cX_b}$ imply $x'x\in{_cX_a}$.

Both the ideals and the subringoids of $X$ are partially
ordered in an obvious manner, and form complete lattices.

\section{Modules and Pointed Overalgebras}
\label{S:Modules}

\subsection{Definitions}
Let $A$ be an algebra of some type $\Omega$. We
define an \emph{abelian group $A$-overalgebra}, or
\emph{$A$-module}, to be an $|A|$-tuple of abelian
groups $\{\,{_aM}\,\}_{a\in A}$, provided with, for
each $n$-ary operation symbol $\omega$ and
each $n$-tuple $\mathbf a$ of elements of $A$, a
group homomorphism $\omega^M_{\mathbf
a}:{_{\mathbf a}M}\to
{_{\omega(\mathbf a)}M}$, where
$_{\mathbf a}M$ stands for
${_{a_1}M}\times\ldots\times{_{a_n}M}$.

A
\emph{pointed $A$-overalgebra} is the same as an
$A$-module, but with pointed sets $_aM$ instead of
abelian groups, and pointed set maps
$\omega^M_{\mathbf a}:{_{\mathbf a}M}\to{_{\omega(\mathbf
a)}M}$ rather than abelian group homomorphisms.  We write
$_a*^M$, or simply $_a*$, for the basepoint of $_aM$. (The
basepoint of
$_{\mathbf a}M$ is
$\langle{_{a_1}*},\ldots,{_{a_n}*}\rangle$.)
An \emph{$A$-overalgebra} is again the same thing, but
with sets $_aM$ rather than pointed sets or abelian
groups. Finally, an \emph{$A$-set} is simply an
$|A|$-tuple of sets.

\subsection{Homomorphisms}
If $S$ and $T$ are $A$-sets, we call an
$|A|$-tuple of functions $_af:{_aS}\to{_aT}$ an
\emph{$A$-function} and write $f:S\to T$.
If $M$ and $M'$ are $A$-modules, a \emph{homomorphism}
$\phi:M\to M'$ is an $A$-function, the components of which
are abelian group homomorphisms, such that
for each
$n$-ary operation symbol $\omega$, each
$n$-tuple $\mathbf a$ of elements of $A$, and each
$\mathbf m\in{_{\mathbf a}M}$, we have
\[
_{\omega(\mathbf a)}\phi(\omega^M_{\mathbf a}(\mathbf m))
=\omega^{M'}_{\mathbf a}
({_{a_1}\phi(m_1)},\ldots,{_{a_n}\phi(m_n)}).
\]
A \emph{homomorphism of pointed $A$-overalgebras} is the
same thing, except that it is an $A$-tuple of pointed set
maps.  A \emph{homomorphism of $A$-overalgebras} is again
the same, but the functions do not have abelian
group structures or basepoints
to preserve.

We denote the category of $A$-modules by $\Ab[A]$,
the category of pointed $A$-overalgebras by $\Pnt[A]$,
and the category of $A$-overalgebras by $\Ov[A]$.

$\Ab[A]$ is an abelian category, if given $\phi$,
$\phi'\in\Ab[A](M,M')$, we define
$_a(\phi+\phi')={_a\phi}+{_a\phi'}$,
$_a(-\phi)=-{_a\phi}$, and
$_a0=0\in\Ab({_aM},{_aM'})$.  See
\S\ref{S:SubsandQuotients} for descriptions of the kernel
and cokernel of a homomorphism.

\begin{example} \label{Ex:Sil}
Let
$B$ be an algebra of the same type as $A$, and let $\pi:B
\to A$ and $\iota:A\to B$ be homomorphisms such that
$\pi\iota=1_A$.
Then we define the pointed $A$-overalgebra $\lsil
B,\pi,\iota\rsil$ by $_a\lsil
B,\pi,\iota\rsil=\pi^{-1}(a)$, $_a*=\iota(a)$, and
$\omega^{\lsil B,\pi,\iota\rsil}_{\mathbf a}(\mathbf
b)=\omega^B(\mathbf b)$.
\end{example}

\begin{example} \label{Ex:P}
Let $\beta\in\Con A$. We
define $\beta^*=\lsil A(\beta),\pi,\iota\rsil$, where
$A(\beta)$ is the subalgebra of $A^2$ given by pairs of
elements related by $\beta$, $\pi:\pair a{a'}\mapsto a$,
and $\iota:a\mapsto\pair aa$. More generally, if
$\alpha$,
$\beta\in\Con A$ are such that $\alpha\leq\beta$, we
define $\mathrm P[\alpha,\beta]$ by $_a\mathrm
P[\alpha,\beta]=\{\,a'/\alpha\mid a\mathrel\beta a'\,\}$,
$_a*^{\mathrm P[\alpha,\beta]}=a/\alpha$, and
$\omega_{\mathbf a}^{\mathrm P[\alpha,\beta]}(\mathbf
x)=\omega^{A/\alpha}(\mathbf x)$.
\end{example}

\subsection{Beck modules}
$\Ab[A]$ is equivalent to the
category of \emph{Beck modules} over $A$ \cite{B}. The
reason for this equivalence is that $\Ab[A]$ is the
category of abelian group objects of $\Ov[A]$, the
category of Beck modules is the category of
abelian group objects of $(\Omega$-$\Alg\downarrow A)$,
the category of algebras (of type
$\Omega$) over
$A$, and $\Ov[A]$ is equivalent
to $(\Omega$-$\Alg\downarrow A)$. (Note that the term
\emph{Beck module} also encompasses abelian group objects
in the category $(\mathbf V\downarrow A)$ where $\mathbf
V$ is a variety of algebras; we will define an equivalent
category $\Ov[A,\mathbf V]$ based on $\Ov[A]$ in the next
section.)

There are two basic
reasons for working in categories derived from $\Ov[A]$
rather than in those derived from
$(\Omega$-$\Alg\downarrow A)$. One reason is that it is
very useful to be able to define
$A$-modules and similar objects $M$ such that the $_aM$
are not necessarily disjoint. The other reason is that if
we were to define, for example, an $A$-module $M$, in
terms of an algebra
$B$ with a homomorphism
$f:B\to A$, we would still need to talk about the abelian
groups $f^{-1}(a)$ for $a\in A$, which we of course
denote by $_aM$. And, indeed, the abelian group
structure (or, pointed set structure, if we talk
about pointed overalgebras) is central to the theory. So
it seems that we would need to talk about the $_aM$
in any case.  Our definition and notation focus attention
where we argue it should be.

In any case, results concerning Beck modules certainly
apply to the category of $A$-modules as we define it.

We mention also that pointed overalgebras have been
studied before \cite{Oz72a, Oz72b, JDHSm}, in the form of
pointed set objects of
$(\mathbf V\downarrow A)$.

\section{The Total Algebra Construction; $\Ab[A,\mathbf
V]$ and $\Pnt[A,\mathbf V]$} \label{S:Total}

Let $M$ be an $A$-module, pointed
$A$-overalgebra, or $A$-overalgebra. We define the
\emph{total algebra} of
$A$ (called by some authors the \emph{semidirect product})
to be the set $A\semitimes M=\{\,\pair am\mid
m\in{_aM}\,\}$, provided with operations
defined by $\omega^{A\semitimes
M}(\pair{a_1}{m_1},\ldots,\pair{a_n}{m_n})=\pair{\omega(
\mathbf a)}{\omega^M_{\mathbf a}(\mathbf m)}$ for each
$n$-ary operation symbol $\omega$.

Associated with the total algebra is a homomorphism
$\pi_M:A\semitimes M\to A$, defined by $\pi_M:\pair
am\mapsto a$.  In the case of an $A$-module or pointed
$A$-overalgebra, there is also a homomorphism
$\iota_M:A\to A\semitimes M$, defined by
$a\mapsto\pair a{{_a0}}$ (or, for
$M$ a pointed $A$-overalgebra, by $a\mapsto\pair
a{{_a*}}$.) Note that $\pi_M\iota_M=1_A$, and that
$\pi_M$ is onto in these cases.  

In the case of an $A$-set $S$, we can still define
$A\semitimes S$ and $\pi_S$, but they are only a set and
a function, rather than an algebra and a homomorphism.

\subsection{Modules, pointed overalgebras, and
overalgebras totally in a variety $\mathbf V$}
 Let $A$ belong to a
variety
$\mathbf V$. Then we say that an $A$-module, pointed
$A$-overalgebra, or
$A$-overalgebra
$M$ is \emph{totally in $\mathbf V$} if
$A\semitimes M$ belongs to $\mathbf V$. We denote
the full subcategory of $\Ab[A]$ of $A$-modules
(pointed $A$-overalgebras, $A$-overalgebras) totally in
$\mathbf V$ by
$\Ab[A,\mathbf V]$ (respectively, by $\Pnt[A,\mathbf V]$,
$\Ov[A,\mathbf V]$).

\begin{remark}
In example~\ref{Ex:Sil}, $\lsil B,\pi,\iota\rsil\in\mathbf
V$ if
$B\in\mathbf V$. In example~\ref{Ex:P}, if $\alpha$,
$\beta\in\Con A$ for $A\in\mathbf V$, then
$\beta^*\in\mathbf V$ and
$\mathrm P[\alpha,\beta]\in\mathbf V$.
\end{remark}

\begin{example}[Free pointed overalgebras] \label{Ex:Free}
There is an evident forgetful functor
\[\mathcal U_A:\Pnt[A,\mathbf V]\to A\text{-}\Set.\]
If
$S$ is an $A$-set, then let $B=A\coprod\mathcal
F_{\mathbf V}(A\semitimes S)$, and
$P=\lsil B,\pi,\iota\rsil$, where $\iota:A\to B$ is the
insertion into the coproduct, and
$\pi:B\to A$ is the
homomorphism determined by $1_A$, $\pi_S:A\semitimes S\to
|A|$, and the universal property of the relatively free
algebra. Then
$\pair P\xi$ is a universal arrow from
$S$ to
$\mathcal U_A$, where $\xi:S\to\mathcal U_AP$ is defined
by $_a\xi:s\mapsto\pair as\in{_aP}$.
We say that $P$ is \emph{free on $S$ (relative to $\mathbf
V$)}.
\end{example}

\section{$A$-operations}\label{S:A-Operations}

If $w$ is an $n$-ary operation on a set $A$, and $M$ is an
$A$-tuple of abelian groups, pointed sets, or simply sets,
we call an
$A^n$-tuple $W$ of abelian group or pointed
set homomorphisms (or, in the case of an $A$-tuple of
sets, functions)
$W_{\mathbf a}:{_{\mathbf a}M}\to{_{w(\mathbf a)}M}$ an
\emph{$A$-operation on $M$, over $w$}.  If we denote
the disjoint union of the $_aM$ by $A\semitimes M$, then
besides possibly consisting of homomorphisms of abelian
groups or of pointed sets, $W$ consists of exactly the
information needed to specify an
$n$-ary operation, which we denote by $A\semitimes W$, on
$A\semitimes M$, such that the obvious function
$\pi:A\semitimes M\to A$ is a homomorphism from the
algebra
$\pair{A\semitimes M}{A\semitimes W}$ to the algebra
$\pair Aw$.

We define $\Clo_nM$ to the the set of pairs $\pair wW$
such that $w$ is an $n$-ary operation on $A$, and $W$ is
an
$A$-operation on
$M$ over
$w$.
We set $\pi^{\Clo M}_{in}=\pair{\pi^{\Clo
A}_{in}}{\tilde\pi_{in}}$, where
$(\tilde\pi_{in})_{\mathbf a}(\mathbf m)=m_i$, and given
$\pair{w_1}{W_1}$, $\ldots$,
$\pair{w_{n'}}{W_{n'}}\in\Clo_nM$ and
$\pair{w'}{W'}\in\Clo_{n'}M$, we define the composite
to be $\pair{w'\mathbf w}{W'\circ_{\mathbf w}\mathbf W}$,
where
$W'\circ_{\mathbf w}\mathbf W$ is the $A$-operation on
$M$ over $w'\mathbf w$ such that for each
$n$-tuple
$\mathbf a$ of elements of $A$, and $n$-tuple $\mathbf m$
of elements of $_{\mathbf a}M$,
\[
(W'\circ_{\mathbf w}\mathbf W)_{\mathbf a}(\mathbf
m)=W'_{\langle w_1(\mathbf a),\ldots,w_{n'}(\mathbf
a)\rangle}((W_1)_{\mathbf a}(\mathbf
m),\ldots,(W_{n'})_{\mathbf a}(\mathbf m)).
\]

\begin{notation} In what folows, we will abbreviate
$\langle w_1(\mathbf a),\ldots,w_{n'}(\mathbf a)\rangle$
by $\mathbf w(\mathbf a)$.
\end{notation}

\begin{theorem} \label{T:IsClone}
We have
\begin{enumerate}
\item $\Clo M$ is a clone; and
\item the $\N$-tuple
of functions defined by $\pair wW\mapsto w$ form a clone
homomorphism $\Lambda_M:\Clo M\to\Clo A$.
\end{enumerate}
\end{theorem}

\begin{proof} 
The only nonobvious part is to show that if
$\pair{w'}{W'}\in\Clo_{n'}M$, and $\pair{w_1}{W_1}$,
$\ldots$, $\pair{w_{n'}}{W_{n'}}\in\Clo_nM$, then
$\pair{w'\mathbf w}{W'\circ_{\mathbf w}\mathbf
W}\in\Clo_nM$, which, in case
$M$ is an $A$-tuple of abelian groups or pointed sets,
requires showing that for each $\mathbf a\in A^n$,
$(W'\circ_{\mathbf w}\mathbf W)_{\mathbf a}$ is a
homomorphism. However, $(W'\circ_{\mathbf w}\mathbf
W)_{\mathbf a}$ is the composition of two homomorphisms:
$W'_{\mathbf w(\mathbf a)}$, and the homomorphism from
$_{\mathbf a}M$ to the product $_{\mathbf w(\mathbf a)}M$
determined by the
$W_{\mathbf a}:{_{\mathbf a}M}\to{_{w_i(\mathbf a)}M}$.
\end{proof}

We call $\Clo M$ the \emph{clone
of
$A$-operations on $M$}.

\subsection{Preservation of $A$-operations}
Let $A$ be a set, and let $M$, $M'$ be $A$-tuples of
abelian groups, pointed sets, or simply sets. Let
$w\in\Clo_nA$ and let $W$, $\bar W$ be such that $\pair
wW\in\Clo_nM$ and $\pair w{\bar W}\in\Clo_nM'$. If $\phi$
is an $A$-tuple of abelian group homomorphisms (pointed
set homomorphisms, functions)
$_a\phi:{_aM}\to{_aM'}$, then we say that \emph{$\phi$
sends $W$ to $\bar W$} if for all $\mathbf a\in A^n$ and
all
$\mathbf m\in{_{\mathbf a}M}$ we have
\[
_{w(\mathbf a)}\phi(W_{\mathbf a}(\mathbf m))
=\bar W_{\mathbf a}
({_{a_1}\phi(m_1)},
\ldots,{_{a_n}\phi(m_n)}).
\]
Of course, we have seen a similar equation before in the
definition of homomorphisms $\phi:M\to M'$ of $A$-modules,
pointed
$A$-overalgebras, and $A$-overalgebras, when $A$ is an
algebra; in such cases, $\phi$ is a homomorphism iff
$\phi$ sends $\omega^M$ to $\omega^{M'}$ for each
$\omega$.

\begin{theorem}\label{T:Preservation} Let $A$ be a set,
$M$,
$M'$
$A$-tuples of abelian groups, pointed sets, or simply
sets,
$w'\in\Clo_{n'}A$, $w_1$, $\ldots$, $w_{n'}\in\Clo_n A$,
and $W'$, $\bar W'$, $W_1$, $\ldots$, $W_{n'}$, $\bar
W_1$, $\ldots$, $\bar W_{n'}$ be such that
$\pair{w'}{W'}\in\Clo_{n'}M$, $\pair{w'}{\bar
W'}\in\Clo_{n'}M'$, and $\pair{w_i}{W_i}\in\Clo_nM$ and
$\pair{w_i}{\bar W_i}\in\Clo_nM'$ for all $i$. Let $\phi$
be an $A$-tuple of abelian group homomorphisms, pointed
set homomorphisms, or functions, respectively. If $\phi$
sends
$W'$ to $\bar W'$ and $W_i$ to $\bar W_i$ for all $i$,
then $\phi$ sends
$W'\circ_{\mathbf w}\mathbf W$ to
$\bar W'\circ_{\mathbf
w}\bar\mathbf W$.
\end{theorem}

\begin{proof} For each $\mathbf a\in A^n$, and each
$\mathbf m\in{_{\mathbf a}M}$, we have
\begin{align*}
_{w'\mathbf w(\mathbf a)}\phi((W'\circ_{\mathbf w}\mathbf
W)_{\mathbf a}(\mathbf m))
&={_{w'\mathbf w(\mathbf a)}\phi(W'_{\mathbf w(\mathbf a)}
((W_1)_{\mathbf a}(\mathbf m),\ldots,(W_{n'})_{\mathbf
a}(\mathbf m)))}\\
&=\bar W'_{\mathbf w(\mathbf a)}(\ldots,{_{w_i(\mathbf a)}
\phi((W_i)_{\mathbf a}(\mathbf m))},\ldots)\\
&=\bar W'_{\mathbf w(\mathbf a)}(\ldots,(\bar
W_i)_{\mathbf a}({_{a_1}\phi(m_1)},\ldots,
{_{a_n}\phi(m_n)}),\ldots)\\
&=(\bar W'\circ_{\mathbf w}\bar\mathbf W)_{\mathbf a}
({_{a_1}\phi(m_1)},\ldots,
{_{a_n}\phi(m_n)}).
\end{align*}
\end{proof}

\subsection{The $A$-operations $t^M$}
Let $A$ be an algebra, and let $M$ be an $A$-module,
pointed
$A$-overalgebra, or $A$-overalgebra.
As we mentioned before, $\Clo\Omega$ is a free clone on
the type $\Omega$. We have an $\N$-tuple of mappings
sending $n$-ary operation symbols $\omega$ to
$\pair{\omega^A}{\omega^M}\in\Clo M$. By the universal
property of $\Clo\Omega$, there is a corresponding clone
homomorphism
$\Phi^M:\Clo\Omega\to\Clo M$. Thus, for each $n$-ary term
$t$, there is an $n$-ary $A$-operation $t^M$ over
$t^A$ such that $\Phi^M_n(t)=\pair{t^A}{t^M}$.

\begin{theorem} Let $A$ be an algebra of type $\Omega$,
and $M$ an $A$-module, pointed $A$-overalgebra, or
$A$-overalgebra.
We have
\begin{enumerate}
\item $\Lambda_M\Phi^M=\Phi^A$;
\item $(\omega(\mathbf x))^M=\omega^M$;
\item for each $n$-ary term $t$, and $n$-tuple
$\pair{a_1}{m_1}$, $\ldots$, $\pair{a_n}{m_n}$ of
elements of $A\semitimes M$, we have
$t^{A\semitimes
M}(\pair{a_1}{m_1},\ldots,\pair{a_n}{m_n})=\pair{t^A(\mathbf
a)}{t^M(\mathbf m)}$;
\item if $M'$ is another $A$-module, pointed
$A$-overalgebra, or $A$-overalgebra, $\phi:M\to M'$
is a homomorphism, and $t\in\Clo_n\Omega$, then $\phi$
sends $t^M$ to $t^{M'}$; and
\item if $B$ is an algebra of the same type as $A$,
and $\pi:B \to A$ and $\iota:A\to B$ are homomorphisms
such that $\pi\iota=1_A$, then for each $n$-tuple of
elements of $A$, and each $\mathbf b\in{_{\mathbf a}\lsil
B,\pi,\iota\rsil}$, we have
$t^{\lsil
B,\pi,\iota\rsil}_{\mathbf a}(\mathbf b)=t^B(\mathbf b)$.
\end{enumerate}
\end{theorem}

\subsection{$A$-operations, identities, and varieties}

\begin{theorem}
Let $M$ be an $A$-module, pointed
$A$-overalgebra, or
$A$-overalgebra. We have
\begin{enumerate}
\item If $t$ and $t'$ are $n$-ary
terms such that
$A$ satisfies the identity $t=t'$, then
$t^M=(t')^M$ iff $t^{A\semitimes M}=(t')^{A\semitimes
M}$; and
\item $M$ is totally in a variety $\mathbf V$ iff
$\Phi^M=\Phi'\Phi_{\mathbf V}$ for some clone
homomorphism $\Phi':\Clo\mathbf V\to\Clo M$, i.e. iff
$\Phi^M$ factors through $\Phi_{\mathbf
V}:\Clo\Omega\to\Clo\mathbf V$.
\end{enumerate}
\end{theorem}

In case the equivalent conditions of part (1) of the
theorem are satisfied, we say that $M$ \emph{satisfies the
identity}
$t=t'$.

\section{Subobjects and
Quotients}\label{S:SubsandQuotients}

\subsection{Factorization of homomorphisms}
If $\phi:M\to M'$ is a homomorphism of $A$-modules, we can
construct a new
$A$-module
$\MyIm\phi$, the
\emph{image} of
$\phi$, by defining $_a\MyIm\phi=\MyIm{_a\phi}$ for all
$a$ and letting the $A$-operations be the restrictions of
those of
$M'$.
$\phi$ then factors in an obvious way as
$\phi_m\phi_e$, where
$\phi_e:M\to\MyIm\phi$ and $\phi_m:\MyIm\phi\to M'$ are
$A$-module homomorphisms such that
for each
$a$,
$_a\phi_e$ is onto and $_a\phi_m$ is one-one. We say that
$\phi_e$ is \emph{onto} and that $\phi_m$ is
\emph{one-one}.

The onto homomorphisms of $A$-modules form a
subcategory
$\mathbf E$ of
$\Ab[A]$, and the one-one
homomorphisms form a subcategory $\mathbf M$.
The pair
$\pair{\mathbf E}{\mathbf M}$ is an example of
a factorization system in $\Ab[A]$. (See \cite{A-H-S} for
a definition of this concept.)

Similar definitions and remarks apply to homomorphisms of
pointed
$A$-overalgebras, $A$-overalgebras, and $A$-sets.

\subsection{Subobjects}
If $M$, $M'$ are $A$-modules (pointed
$A$-overalgebras, $A$-overalgebras, $A$-sets) then we
say that $M'\leq M$, or that $M'$ is a
\emph{submodule} (respectively, \emph{sub pointed
overalgebra}, \emph{sub overalgebra}, \emph{$A$-subset})
of $M$, if
$_aM'\subseteq{_aM}$ for every $a$. For a given $M$,
the subobjects form a complete lattice, which we denote by
$\Sub M$.

On the other
hand, given two one-one homomorphisms $\iota$,
$\iota'$ of $A$-modules (or pointed overalgebras, etc.)
with codomain $M$, we say that $\iota\leq\iota'$ if there
is a homomorphism
$\phi$ such that $\iota=\iota'\phi$. This results in a
preorder, and as usual with a preorder, we say that
$\iota$ and $\iota'$ are \emph{equivalent} if
$\iota\leq\iota'$ and $\iota'\leq\iota$. The resulting
partially-ordered set of equivalence classes is
isomorphic to $\Sub M$.

 We will speak of $A$-submodules
of a given module
$M$,
\emph{generated by} an $A$-subset $S$ of $M$. This
means the smallest $A$-submodule $M'$, in the lattice of
submodules of $M$, such that $_aS\subseteq{_aM'}$ for
each $a\in A$.

\begin{theorem} Let $M$ be an $A$-module or pointed
$A$-overalgebra. We have
\begin{enumerate}
\item If $M'\leq M$ is a submodule or sub pointed
overalgebra, respectively, then
$A\semitimes M'$ is a subalgebra of $A\semitimes M$; and
\item 
$\Sub M\cong\I_{\Sub A\semitimes M}[A,\top]$
under the mapping $M'\mapsto A\semitimes M'$.
\end{enumerate}
\end{theorem}

\begin{example} Let $\phi:M\to M'$ be a
homomorphism of $A$-modules. We define
$\Ke\phi$ by
$_a(\Ke\phi)=\Ke({_a\phi})$, and this can easily be shown
to be a submodule of $M$, and to be a kernel of $\phi$
in the sense of the theory of additive categories.
\end{example}

\subsection{Quotient objects}
If $M$ is a pointed $A$-overalgebra or
$A$-overalgebra, a
\emph{congruence} of
$M$ is an
$|A|$-tuple $\gamma$ of equivalence relations $_a\gamma$
such that for each $n$-ary $\omega$, each $n$-tuple
$\mathbf a$ of elements of $A$, and each $\mathbf
m$, $\mathbf m'\in{_{\mathbf a}M}$ such that
$m_i\mathrel{{_{a_i}\gamma}}m'_i$ for all $i$, we have
$\omega^M_{\mathbf a}(\mathbf m)\mathrel{{_{\omega(\mathbf
a)}\gamma}}\omega^M_{\mathbf a}(\mathbf m')$. If $M$ is an
$A$-module, then we impose the
additional condition that each $_a\gamma$ be an abelian
group congruence of $_aM$.
The congruences of $M$ are partially ordered in an
obvious manner, and this gives rise to a complete lattice
of congruences, $\Con M$.

If $M$ and $\gamma\in\Con M$ are given, then there is a
unique structure of $A$-module (pointed $A$-overalgebra,
$A$-overalgebra) on the $|A|$-tuple of abelian groups
(respectively, pointed sets, sets) $_aM/{_a\gamma}$,
which we denote by $M/\gamma$, such that the natural maps
$\nat{_a\gamma}$ form a homomorphism $\nat\gamma:M\to
M/\gamma$ of
$A$-modules (respectively, of pointed $A$-overalgebras,
$A$-overalgebras).

Just as $\Sub M$ is isomorphic to a lattice of
equivalence classes of one-one homomorphisms with
codomain $M$, $\Con M$ is isomorphic to a lattice of
equivalence classes of onto homomorphisms with domain $M$.

\begin{theorem} Let $A$ be an algebra, and let $M$
be an $A$-module, pointed $A$-overalgebra, or
$A$-overalgebra.
\begin{enumerate}
\item If $\gamma\in\Con M$, then the binary relation
$A\semitimes\gamma$, defined by $\pair
am\mathrel{A\semitimes\gamma}\pair {a'}{m'}$ iff $a=a'$
and
$m\mathrel{{_a\gamma}} m'$, is a congruence of
$A\semitimes M$; and
\item
$\Con M\cong\I_{\Con A\semitimes M}[\bot,\ker\pi_M]$, via
the mapping $\gamma\mapsto A\semitimes\gamma$.
\end{enumerate}
\end{theorem}

Just as a congruence on an abelian group is determined by
the equivalence class of the $0$ element, a congruence
$\gamma$ on an $A$-module $M$ is determined by the
submodule $M'$ such that each $_aM'$ is the
$_a\gamma$-class of
$_a0^M$.

\begin{example}
Let $\phi:M\to M'$ be a homomorphism of
$A$-modules. We define $\Co\phi$ by
$_a(\Co\phi)=\Co({_a\phi})$. It is straightforward to
prove that there is a unique structure of $A$-module on
$\Co\phi$, making the natural maps from $M'$ onto
$\Co\phi$ a homomorphism, and that
$\Co\phi$ is a cokernel of
$\phi$.
\end{example}

\section{The Clone of
Polynomials} \label{S:Polynomials}

Let $A$ be an algebra in a variety $\mathbf V$ of algebras
of type $\Omega$. We define the \emph{algebra of $n$-ary
polynomials with coefficients in $A$ (relative to $\mathbf
V$)} to be the  coproduct $A\coprod\mathcal F_{\mathbf
V}\{\,x_1,\ldots, x_n\,\}$. Of course, this
algebra is defined only up to isomorphism,
but we choose one such algebra, and denote it by
$Q_n$. For each $n$, there is a
distinguished one-one homomorphism from $A$ to
$Q_n$, and we will use the same notation
$\iota_A$ for any of these homomorphisms.  For $a\in A$,
we will typically write
$a$ for $\iota_A(a)$ in our formulas.

 If
$\Pi\in Q_n$, and
$\mathbf a$ is an
$n$-tuple of elements of $A$, then we define
$\Pi^A(\mathbf a)$ to be the image of $\Pi$ in $A$, under
the unique homomorphism
$\phi:Q_n\to A$ sending each $a\in A$ to
$a$ and each $x_i$ to $a_i$. Such a homomorphism exists,
and is unique, by the universal property of
$Q_n$. Similarly, if $B$ is an algebra in
$\mathbf V$ with a distinguished homomorphism $f:A\to B$,
and
$\mathbf b$ is an $n$-tuple of elements of $B$, then we
define $\Pi^{B,f}(\mathbf b)$ to be the image in $B$
of
$\Pi$, under the unique homomorphism
$\phi:Q_n\to B$ sending each $a\in A$ to
$f(a)$, and each $x_i$ to $b_i$.

We will write $x$ for $x_1$ when
talking about unary polynomials, just as we do with unary
terms.

\begin{theorem}\label{T:PolyDef} Let $A$ be an algebra in
a variety
$\mathbf V$ of algebras of type $\Omega$. We have
\begin{enumerate}
\item If we
define $\Pol_n(A,\mathbf V)=|Q_n|$,
$\pi^{\Pol(A,\mathbf V)}_{i,n}=x_i$ and, given
$n$, $n'$, an $n'$-tuple $\pmb\Pi$ of $n$-ary
polynomials, and an $n'$-ary polynomial $\Pi'$,
$\Pi'\pmb\Pi=(\Pi')^{Q_n,\iota_A}(\pmb\Pi)$,
then $\Pol(A,\mathbf V)$ is a clone;
\item The mappings $(\Phi^A_{\Pol})_n:\Pol_n(A,\mathbf V)
\to\Clo_n|A|$, given by $\Pi\mapsto\Pi^A$, form a clone
homomorphism;
\item if $B$ is an algebra, provided with a homomorphism
$f:A\to B$, the mappings $\Phi^{B,f}_n:\Pol_n(A,\mathbf
V)\to\Clo_n|B|$, given by $\Pi\mapsto \Pi^{B,f}$, form
a clone homomorphism; and
\item there is a natural clone homomorphism
$\tilde\Phi:\Clo\Omega\to\Pol(A,\mathbf V)$, such that
$\Phi^{B,f}\tilde\Phi=\Phi^B$.
\end{enumerate}
\end{theorem}

\begin{proof} (1), (2), and (3) are
straightforward. For (4), $\tilde\Phi_n$ is simply the
insertion of the relatively free algebra into the
coproduct.
\end{proof}

\begin{remark} $\tilde\Phi_n$ may not be one-one.
\end{remark}

\subsection{Polynomials and polynomial functions} 
We have defined a clone $\Pol(A,\mathbf V)$ of
polynomials with coefficients in $A$ (relative to
$\mathbf V$), and, given an algebra $B\in\mathbf V$ with
a homomorphism $f:A\to B$, a clone homomorphism
$\Phi^{B,f}:\Pol(A,\mathbf V)\to\Clo|B|$, which sends each
$n$-ary polynomial $\Pi$ to $\Pi^{B,f}$. The image
of this clone homomorphism is a subclone of $\Clo|B|$,
the \emph{clone of polynomial functions in $B$ with
coefficients in $A$}.
As with terms and term functions, it is necessary to
distinguish between these two clones.

\subsection{Algebras of polynomials and free pointed
overalgebras}
Recall (example~\ref{Ex:Free}) our recipe for a pointed
overalgebra free on an
$A$-set. An $n$-tuple $\mathbf a$ of elements of $A$
determines an
$A$-set $S$ with disjoint union
$\{\,x_1,\ldots,x_n\,\}$, by the rule
$x_i\in{_{\bar a}S}\iff a_i=\bar a$. Then if $\pi_{\mathbf
a}$ is evaluation at $\mathbf a$, i.e., the homomorphism
$\Pi\mapsto\Pi^A(\mathbf a)$, we have
$\lsil\Pol_n(A,\mathbf V),\pi_{\mathbf a},\iota_A\rsil$ a
pointed $A$-overalgebra, totally in $\mathbf V$ and free
on $S$. We will use this fact when we use the algebra of
unary polynomials in our construction of the enveloping
ringoid, and just below in the definition of the
$A$-operation
$\Pi^P$ for $P$ a pointed $A$-overalgebra.

\subsection{Polynomials and pointed overalgebras}
Given
$\Pi\in\Pol_n(A,\mathbf V)$,
$P\in\Pnt[A,\mathbf V ]$, $\mathbf a\in A^n$, and
$\mathbf p\in {_{\mathbf a}P}$, we define
$\Pi^P_{\mathbf a}(\mathbf p)$ to be the image of $\Pi$
in $P$, under the unique homomorphism from
$\lsil\Pol_n(A,\mathbf V),\pi_{\mathbf a},\iota_A\rsil$
to $P$ sending each $x_i$ to $p_i$.

\begin{theorem} \label{T:PolyAction} Let $A$ be an
algebra, and $P\in\Pnt(A,\mathbf V)$. We have
\begin{enumerate}

\item The $\N$-tuple of mappings sending $n$-ary
polynomials $\Pi$ to $\pair{\Pi^A}{\Pi^P}$ form a clone
homomorphism $\Phi_{\Pol}^P:\Pol(A,\mathbf V)\to\Clo M$,
such that $\Lambda_P\Phi_{\Pol}^P=\Phi^A_{\Pol}$ and
$\Phi^P_{\Pol}\tilde\Phi=\Phi^P$;

\item for each $n$-ary polynomial $\Pi$, and $n$-tuple
$\pair{a_1}{p_1}$, $\ldots$, $\pair{a_n}{p_n}\in
A\semitimes P$, we have
$\Pi^{\pair{A\semitimes
P}{\iota_P}}(\pair{a_1}{p_1},\ldots,\pair{a_n}{p_n})=
\pair{\Pi^A(\mathbf a)}{\Pi^P_{\mathbf a}(\mathbf p)}$;

\item if $P'$ is another pointed $A$-overalgebra,
$\phi:P\to P'$ a homomorphism, and $\Pi\in\Pol_n(A,\mathbf
V)$, then $\phi$ sends $\Pi^P$ to
$\Pi^{P'}$; and

\item if $B \in\mathbf V$, $\pi:B\to A$, and
$\iota:A\to B$ are such that $\pi\iota=1_A$,
then $\Pi^{\lsil B,\pi,\iota\rsil}_{\mathbf a}(\mathbf b)
=\Pi^{B,\iota}(\mathbf b)$.

\end{enumerate}
\end{theorem}

\begin{proof}
All parts are straightforward, if not
obvious, except perhaps (3), the key to proving
which is to know that if
$\omega\in\Clo_{n'}\Omega$, then $\phi$ sends
$\tilde\Phi(\omega)^P=\omega^P$ to
$\omega^{P'}=\tilde\Phi(\omega)^{P'}$ because it is a
homomorphism, and then to use
theorem~\ref{T:Preservation}.\end{proof}

\subsection{Polynomials and $A$-modules}
There is a natural forgetful functor $\mathcal
U:\Ab[A,\mathbf V]\to\Pnt[A,\mathbf V]$, which takes an
$A$-module $M$ to the pointed $A$-overalgebra $\mathcal
UM$ defined by $_a\mathcal UM={_aM}$, $_a*^{\mathcal
UM}={_a0^M}$. If $\Pi\in\Pol_n(A,\mathbf V)$ and $\mathbf
a\in A^n$, we define $\Pi^M_{\mathbf a}=\Pi^{\mathcal
UM}_{\mathbf a}$.

\begin{theorem} Let $M$ be an $A$-module totally in
$\mathbf V$, $\Pi\in\Pol_n(A,\mathbf V)$, and $\mathbf
a\in A^n$. Then $\Pi^M_{\mathbf a}$ is a group
homomorphism and $\pair{\Pi^A}{\Pi^M}\in\Clo_nM$.
\end{theorem}

\begin{proof} We
use the fact that $\Pol_n(A,\mathbf V)$ is generated by
the projections and the constants $a\in A$. Certainly the
conclusion is true for $\Pi$ one of these generators.
Suppose it is true for the
components of $\pmb\Pi\in\Pol_n(A,\mathbf V)^{n'}$, and
let
$\omega\in\Omega_{n'}$. We have
$(\Phi^M_{\Pol})_{n'}(\omega)=\pair{\omega^A}{\omega^M}
\in\Clo_{n'}M$, and
$\pair{\Pi^A_i}{\Pi^M_i}\in\Clo_nM$ for all $i$. Thus, by
theorem~\ref{T:IsClone},
$\tilde\Phi
(\omega)^M\circ_{\pmb{\Pi}^A}\pmb\Pi^M$ is an
$A$-operation on $M$ over $\tilde\Phi(\omega)\pmb\Pi$.
But, we have
\begin{align*}
(\tilde\Phi(\omega)^M\circ_{\pmb\Pi^A}\pmb\Pi^M)_{\mathbf
a}(\mathbf m)
&=(\tilde\Phi(\omega)^{\mathcal
UM}\circ_{\pmb
\Pi^A}\pmb\Pi^{\mathcal UM})_{\mathbf a}(\mathbf m)\\
&=(\tilde\Phi(\omega)\pmb\Pi)^{\mathcal
UM}_{\mathbf a}(\mathbf m)\\
&=(\tilde\Phi(\omega)\pmb\Pi)^M_{\mathbf a}(\mathbf m)
\end{align*}
for every $\mathbf m\in{_{\mathbf a}M}$, by definition
and by the fact that $\Phi_{\Pol}^{\mathcal UM}$ is a
clone homomorphism. Thus, the conclusion is true for
$\Pi=\tilde\Phi(\omega)\pmb\Pi$.
\end{proof}.

\section{$Z_M$}\label{S:ZsubM}
Let $A$ be a set, and $M$ an $A$-tuple of abelian
groups. If $w$ is an $n$-ary operation on $A$,
$W$ is an $n$-ary $A$-operation on $M$ over $w$,
$\mathbf a$ is an $n$-tuple of elements of $A$, and
$1\leq i\leq n$, then we denote by
$W_{\mathbf a,i}$ the abelian group
homomorphism from $_{a_i}M$ to $_{w(\mathbf
a)}M$ defined by $m\mapsto W_{\mathbf
a}({_{a_1}0},\ldots,{_{a_{i-1}}0},m,{_{a_{i+1}}0},
\ldots,
{_{a_n}0})$.

\begin{theorem} Let $A$ be a set with an $n$-ary
operation $w$, and $M$ and $M'$ $A$-tuples of abelian
groups. Let $W$ be an $A$-operation over $w$ on
$M$, and $W'$ an $A$-operation over $w$ on $M'$.
\begin{enumerate}
\item If $M=M'$, then $W=W'$ iff $W_{\mathbf
a,i}=W'_{\mathbf a,i}$ for all $\mathbf a$ and $i$.
\item If $\phi$ is an $A$-tuple of homomorphisms
$_a\phi:{_aM}\to{_aM'}$, then $\phi$ sends $W$ to $W'$
iff for all $i$ such that $1\leq i\leq n$, $\mathbf a\in
A^n$, and $m\in{_{a_i}M}$, we have
\[
\phi_{w(\mathbf a)}(W_{\mathbf a,i}(m))=W'_{\mathbf
a,i}(\phi_{a_i}(m)).
\]
\end{enumerate}
\end{theorem}

If $A$ is an algebra and $M$ is an $A$-module, then the
$\omega^M_{\mathbf a,i}$ are of particular interest.
We denote by $Z_M$ the smallest subringoid of
$\End(M)$ containing all of the $\omega^M_{\mathbf
a,i}$.

\begin{theorem} Let $M$ be an $A$-module. We have
\begin{enumerate}
\item In order for
an $A$-tuple of subgroups $_aM'\subseteq{_aM}$ to be a
submodule, it is necessary and sufficient that
for
each $n$-ary operation symbol $\omega$, each $n$-tuple
$\mathbf a$ of elements of $A$, and each $i$ such that
$1\leq i\leq n$, we have
$\omega_{\mathbf
a,i}({_{a_i}M'})\subseteq{_{\omega(\mathbf a)}M'}$; and
\item
if $M'$ is the submodule of $M$ generated by an
$A$-subset
$S$, then for each $a\in A$, $_aM'$ is the
subgroup of $_aM$ generated by elements of the form $rm$,
where $r\in{_a{(Z_M)}_b}$ for some $b\in A$, and
$m\in{_bS}$.
\end{enumerate}
\end{theorem}

\section{Modulization of a Pointed Overalgebra}
\label{S:Modulization}

Recall (\S\ref{S:Polynomials}) that we denote by $\mathcal
U$ the natural forgetful functor from $\Ab[A]$ to
$\Pnt[A]$, which assigns to each
$A$-module
$M$ the pointed overalgebra $\mathcal UM$ such that
$_a(\mathcal UM)={_aM}$ and $_a*^{\mathcal
UM}={_a0^M}$. For arrows, $\mathcal UM$ sends a
homomorphism $\phi:M\to M'$, consisting of an $A$-tuple of
abelian group homomorphisms $_a\phi:{_aM}\to{_aM'}$, to
the same
$A$-tuple of functions, which are also pointed set
homomorphisms and comprise a pointed overalgebra
homomorphism.
 Note that if $\mathbf V$ is a variety of
algebras, $\mathcal U$ sends
$\Ab[A,\mathbf V]$ into
$\Pnt[A,\mathbf V]$ and we can view $\mathcal U$ as a
functor from $\Ab[A,\mathbf V]$ to $\Pnt[A,\mathbf V]$.

In this section, we will construct universal arrows to all
of these functors. In \S\ref{S:Enveloping}, we will use
the resulting left adjoint functors to construct
enveloping ringoids.

We will call the left adjoint to $\mathcal
U:\Ab[A]\to\Pnt[A]$ the functor of \emph{modulization}
and denote it by $\mathcal M$.
However, we will begin by defining an easier functor
$\hat\mathcal M:\Pnt[A]\to\Ab[A]$. If $P$ is a pointed
$A$-overalgebra, then for each $a\in A$, $_a(\hat\mathcal
MP)$ will be the free abelian group on
$_aP$. We will identify each $p\in{_aP}$ with the
corresponding generator of $_a(\hat\mathcal MP)$. If
$\omega$ is an
$n$-ary operation symbol, and
$\mathbf a$ an $n$-tuple of elements of $A$, then we
will define each $\omega^{\hat\mathcal M(P)}_{\mathbf
a,i}$ on the generators of $_{a_i}\hat\mathcal M(P)$ by
the equation
\[
\omega^{\hat\mathcal MP}_{\mathbf a,i}(p)
=\omega^P({_{a_1}*},\ldots,{_{a_{i-1}}*},p,
{_{a_{i+1}}*},\ldots,
{_{a_n}*});
\]
the $\omega^{\hat\mathcal MP}_{\mathbf a,i}$ then
determine the $A$-operations $\omega^{\hat\mathcal MP}$.

In order to obtain the left adjoint functor $\mathcal M$,
we must have a unit natural transformation
$\eta:1_{\Pnt[A]}\natur\mathcal U\mathcal M$. We have a
natural definition for functions
$_a(\hat\eta P):{_aP}\to{_a(\mathcal U\hat\mathcal MP)}$;
we simply send each $p\in{_aP}$ to itself, viewed as an
element of $_a(\hat\mathcal MP)$. However, the
$A$-tuple of functions $\hat\eta P$ is not a homomorphism
of pointed $A$-overalgebras, because it does not send
each $_a*$ to $_a0$ and does not satisfy the equations
\[
_{\omega(\mathbf a)}(\hat\eta P)(\omega^P_{\mathbf
a}(\mathbf p))
=\omega^{\mathcal U\hat\mathcal MP}_{\mathbf
a}({_{a_1}}(\hat\eta P)(p_1),\ldots,{_{a_n}}(\hat\eta
P)(p_n))
\]
for each $n$-ary operation symbol $\omega$, each
$n$-tuple $\mathbf a$ of elements of $A$, and
each $\mathbf p\in{_{\mathbf a}P}$.

Accordingly, we define $\mathcal MP$ to be the
quotient of $\hat\mathcal MP$ by the smallest
submodule $\mathcal KP$ such that the composite
$A$-function
$\eta P$ of the
$A$-function $\hat\eta P$ and the natural homomorphism
$\nu:\hat\mathcal MP\to\mathcal MP$ yields
a homomorphism of pointed
$A$-overalgebras.  Let $\mathcal SP$ be the smallest
$A$-subset of $\hat\mathcal MP$ such that each
$_a\mathcal SP$ contains
$_a*$, and such that for each $n$-ary operation
symbol $\omega$, each $n$-tuple $\mathbf a$ of
elements of
$A$, and each $\mathbf
p\in{_{\mathbf a}P}$,
$_{\omega(\mathbf a)}\mathcal SP$ contains
$\omega^P_{\mathbf a}(\mathbf
p)-\sum_{i=1}^n\omega^{\hat\mathcal MP}_{\mathbf a
,i}(p_i)=\omega_{\mathbf a}^P(\mathbf
p)-\omega^{\hat\mathcal MP}_{\mathbf a}(\mathbf p)$. Then
$\mathcal KP$ is the submodule of $\hat\mathcal MP$
generated by the
$A$-set
$\mathcal SP$.

\begin{theorem} For all $P\in\Pnt[A]$, $\eta
P$ is a homomorphism of pointed
$A$-overalgebras, and $\pair{\mathcal MP}{\eta P}$ is a
universal arrow from $P$ to
$\mathcal U$.
\end{theorem}

\begin{proof} Let $\omega\in\Omega_n$, $\mathbf a\in
A^n$, and $\mathbf p\in{_{\mathbf a}P}$. For each $a\in
A$, $_a\eta P$ preserves the basepoint, because
$_a*^P\in{_a\mathcal SP}$, so that
\[
_a\eta P({_a*^P})=\nu({_a*^P})={_a0^{\mathcal MP}}.
\]
Also, because $\omega^P_{\mathbf a}(\mathbf
p)-\omega^{\mathcal MP}_{\mathbf a}\in{_{\omega(\mathbf
a)}\mathcal SP}$, we have
\begin{align*}
\omega^{\mathcal U\mathcal MP}_{\mathbf a}({_{a_1}\eta
P(p_1)},\ldots,{_{a_n}\eta P(p_n)})
&=\omega_{\mathbf a}^{\mathcal
MP}({_{a_1}\nu(p_1)},\ldots,{_{a_n}\nu(p_n)})\\
&={_{\omega(\mathbf a)}\nu\left(\omega_{\mathbf
a}^{\hat\mathcal MP}(\mathbf p)\right)}\\
&={_{\omega(\mathbf a)}\nu\left(\omega^P_{\mathbf
a}(\mathbf p)\right)}
\\
&=({_{\omega(\mathbf a)}\nu})({_{\omega(\mathbf
a)}\hat\eta P})(\omega^P_{\mathbf a}(\mathbf p))\\
&={_{\omega(\mathbf a)}\eta P(\omega^P_{\mathbf
a}(\mathbf p))}.
\end{align*}
Thus, $\eta P$ is a homomorphism.

Now let $M\in\Ab[A]$, and $\zeta:P\to\mathcal UM$. We must
show that there is a unique $A$-module homomorphism
$\xi:\mathcal MP\to M$ such that $\zeta=(\mathcal
U\xi)(\eta P)$. By the construction of $\hat\mathcal MP$,
there is a unique $A$-module homomorphism
$\hat\xi:\hat\mathcal MP\to M$ such that $\zeta=(\mathcal
U\hat\xi)(\hat\eta P)$; this shows that $\xi$ is unique if
it exists, because then $\xi\nu=\hat\xi$,
and $\nu$ is an epimorphism.

To show $\xi$ exists, we must show that for
each $a$, $_a\hat\xi$ is zero on $_a\mathcal SP$.
But, we have $_a\hat\xi({_a*^P})={_a0^M}$, because
$_a\zeta({_a*^P})={_a*^{\mathcal UM}}={_a0^M}$, and
if $\omega$ is an $n$-ary operation symbol, $\mathbf a\in
A^n$, and
$\mathbf p\in{_{\mathbf a}P}$, we have
\begin{align*}
_{\omega(\mathbf a)}\hat\xi\left(\omega^P_{\mathbf
a}(\mathbf p) -\sum_{i=1}^n\omega_{\mathbf
a,i}^{\hat\mathcal MP}(p_i)\right)
&={_{\omega(\mathbf a)}\hat\xi}(\omega^P_{\mathbf
a}(\mathbf p))
-\sum^n_i{_{\omega(a)}\hat\xi}(\omega^{\hat\mathcal
MP}_{\mathbf a,i}(p_i)) \\
&=\omega^M_{\mathbf
a}({_{\mathbf a}\hat\xi(\mathbf p)})
-\sum_{i=1}^n\omega^M_{\mathbf a,i}({_{a_i}\hat\xi(p_i)})
\\
&=\omega^M_{\mathbf a}({_{\mathbf a}\zeta(\mathbf p)})
-\sum_{i=1}^n\omega^M_{\mathbf a,i}({_{a_i}\zeta(p_i)}) \\
&={_{\omega(\mathbf a)}0}.
\end{align*}
\end{proof}

It is also obvious from the construction that we have

\begin{theorem} \label{T:LinearComb} Let $P$ be a pointed
$A$-overalgebra, and $a\in A$. Then every element of
$_a(\mathcal MP)$ is a
$\Z$-linear combination of images of elements of $_aP$
under $_a(\eta P)$.
\end{theorem}

If $\mathbf V$ is a variety of algebras to which $A$
belongs, we also want a left adjoint
functor to
the restriction of $\mathcal U$ to $\Ab[A,\mathbf V]$. It
suffices to have
$\mathcal MP$ always be totally in $\mathbf V$ if $P$
is; then the restriction of $\mathcal M$ to
$\Pnt[A,\mathbf V]$ will be the desired left adjoint.

\begin{theorem} We have
\begin{enumerate}
\item $\mathcal MP$ satisfies any identity
satisfied by
$P$.
\item
If $P$ is a pointed $A$-overalgebra totally in $\mathbf
V$, then $\mathcal MP$ is totally in $\mathbf V$.
\item If $P$ is totally in $\mathbf V$,
then $\pair{\mathcal MP}{\eta P}$ is a universal arrow
from $P$ to the restriction of $\mathcal U$ to
$\Ab[A,\mathbf V]$.
\end{enumerate}
\end{theorem}

\begin{proof} (3) follows from (2), which follows from
(1). To prove (1), first we prove that for each
$n$-ary term $t$,
$n$-tuple $\mathbf a$ of elements of $A$,
$i$ such that $1\leq i\leq n$, and
$p\in{_{a_i}P}$, we have
\begin{equation} \label{E:Termfact}
t^{\hat\mathcal MP}_{\mathbf a,i}(p)-t^P_{\mathbf a}
({_{a_1}}*,\ldots,{_{a_{i-1}}}*,p,{_{a_{i+1}}}*,\ldots,
{_{a_n}}*)\in{_{t^A(\mathbf a)}\mathcal KP};
\end{equation}
we will use the fact that $\Clo_n\Omega$ is generated by
the $x_i$. Let
$j\in\{\,1,\ldots,n\,\}$. We have

\begin{align*}
(x_j)^{\hat\mathcal MP}_{\mathbf
a,i}(p)
-(x_j)^P_{\mathbf
a}&({_{a_1}*},\ldots,
{_{a_{i-1}}*},p,{_{a_{i+1}}*},\ldots,{_{a_n}*}) \\
&=\begin{cases} p-p,	&\text{if $i=j$, and}\\
-{_{a_j}*},	&\text{if $i\neq j$}
\end{cases}\\
&\in{_{a_j}\mathcal KP}.
\end{align*}

Now, suppose
$\omega$ is an $m$-ary operation symbol, and
\eqref{E:Termfact} is true for
$n$-ary terms $s_1$,
$\ldots$,
$s_m$. Then we will show it is true for the
$n$-ary term $\omega\mathbf s$.
We have
\begin{align*}
(\omega\mathbf s)^{\hat\mathcal MP}_{\mathbf a,i}(p)
&=
\omega^{\hat\mathcal MP}_{\mathbf s(\mathbf a)}
\left((s_1)_{\mathbf a,i}^{\hat\mathcal MP}(p),\ldots,
(s_m)_{\mathbf a,i}^{\hat\mathcal MP}(p)\right) \\
&=\sum_{j=1}^m\omega_{\mathbf s(\mathbf
a),j}^{\hat\mathcal MP} \left((s_j)_{\mathbf a,i}^{\hat
MP}(p)\right) \\
&=\sum_{j=1}^m\omega^{\hat\mathcal MP}_{\mathbf
s(\mathbf a),j} \left((s_j)_{\mathbf
a}^P({_{a_1}*},\ldots,{_{a_{i-1}}*},p,{_{a_{i+1}}*},
\ldots,{_{a_n}*})+k_j\right)
\\
&=\left(\sum_{j=1}^m\omega^{\hat\mathcal MP}_{\mathbf s(
\mathbf a),j}\left((s_j)^P_{\mathbf
a}({_{a_1}*},\ldots,{_{a_{i-1}}*},p,{_{a_{i+1}}*},\ldots,
{_{a_n}*}
)\right)\right) + k \\ 
&=\omega^P_{\mathbf s(\mathbf a)}\left(\ldots,
(s_j)^P_{\mathbf a}({_{a_1}*},\ldots,{_{a_{i-1}}*},p,
{_{a_{i+1}}*},\ldots,{_{a_n}*}),\ldots\right)+k' \\
&=(\omega\mathbf s)^P_{\mathbf
a}({_{a_1}*},\ldots,{_{a_{i-1}}*},p,{_{a_{i+1}}*},\ldots,
{_{a_n}*})
+ k'
\end{align*}
where each $k_j\in{_{s_j(\mathbf a)}\mathcal KP}$, and
$k$,
$k'\in{_{\omega\mathbf s(\mathbf a)}\mathcal KP}$

Now, each $A$-operation $t^{\hat\mathcal MP}$ is
determined by the actions of the
$t^{\hat\mathcal MP}_{\mathbf a,i}$ on elements of
$_{a_i}P$.  If $t$ and $t'$ are $n$-ary terms such that
$t^P=(t')^P$, then $t^A=(t')^A$ because
$\pi_P:A\semitimes P\to P$ is onto, and \eqref{E:Termfact}
implies that values of
$t^{\hat\mathcal MP}$ and
$(t')^{\hat\mathcal MP}$ differ by elements of
the submodule $\mathcal KP$. This implies that
$t^{\mathcal MP}=(t')^{\mathcal MP}$.
In other words,
$\mathcal MP$ satisfies any identity satisfied by $P$.
\end{proof}

\section{Enveloping Ringoids} \label{S:Enveloping}

In this section, we will use the modulization
functor $\mathcal M:\Pnt[A,\mathbf V]\to\Ab[A,\mathbf
V]$ to construct the enveloping ringoid
$\Z[A,\mathbf V]$ of an algebra
$A$ in a variety $\mathbf V$.

\newcommand\myiota{\iota_A}
If $A$ is an algebra in a variety $\mathbf V$, we will
abbreviate
$\Pol_1(A,\mathbf V)$ by $U$, and continue to denote the
insertion of
$A$ into $U$ by $\myiota$. If
$b\in A$, we will write
$U_b$ for the pointed overalgebra $\lsil
U,\pi_b,\myiota\rsil$, where $\pi_b:U\to A$
is the homomorphism given by $u\mapsto u(b)$.
If $a\in A$, we will write $_aU_b$ for $_a(U_b)$,
and if $\mathbf a$ is an
$n$-tuple of elements of
$A$, we will write $_{\mathbf a}U_b$ for
${_{a_1}U_b}\times\ldots\times{_{a_n}U_b}$.

Now, as discussed previously in
\S\ref{S:Polynomials}, $U_b$ is a free pointed
$A$-overalgebra, totally in $\mathbf V$, on one generator,
$x\in{_bU_b}$.  In other words,
$U_b$ is free on the $A$-subset $S$ such that
$_bS=\{\,x\,\}$ and $_aS=\varphi$ for $a\neq b$.
From the fact that $\mathcal M$ is left adjoint to
$\mathcal U$, it follows that
$\mathcal MU_b\in\Ab[A,\mathbf V]$ is free
on $_b(\eta U_b)(x)$, where $\eta$ is the unit
natural transformation, described in
\S\ref{S:Modulization}, of the adjunction between
$\mathcal M$ and $\mathcal U$.

Let $P$ be a pointed $A$-overalgebra totally in $\mathbf
V$. Since $U_b$ is free on $x\in{_bU_b}$, if $p\in{_bP}$,
then there is a unique homomorphism $\psi^P_{b,p}:U_b\to
P$ taking $x$ to $p$. We will denote this
homomorphism by $\psi_{b,p}$ when there is no
confusion. Recall that by the definition in
\S\ref{S:Polynomials}, if $p\in{_bP}$, and $u$
is a unary polynomial relative to $\mathbf V$,
then $_{u(b)}(\psi_{b,p})(u)=u^P_{\langle
b\rangle}(p)$.

Similarly, suppose
$M$ is an
$A$-module totally in $\mathbf V$. Since $\mathcal MU_b$
is free on $_b(\eta U_b)(x)$, if $m\in{_bM}$, there is
a unique homomorphism
$\theta_{b,m}=\theta^M_{b,m}:\mathcal M U_b\to M$ taking
$_b(\eta U_b)(x)$ to $m$. We have
$(\theta^M_{b,m})(\eta U_b)=\psi^{\mathcal UM}_{b,m}$.

We define $_aZ_b={_a(\mathcal MU_b)}$ for $a$, $b\in A$.
If $u\in U$ and $b\in A$, we
also define
$(u)_b$ to mean
${_{u(b)}(\eta U_b)(u)}$. That is, $(u)_b$ is the
equivalence class of $u$ in the modulization of $U_b$.
 Let
$M$ be an $A$-module totally in $\mathbf V$. If
$z\in{_aZ_b}$, and $m\in{_bM}$, then we define $zm$ to be
$_a(\theta_{b,m})(z)$.
Given $z\in{_aZ_b}$ and $z'\in{_bZ_c}$, we define $zz'$
to be the result of the action of $z$ on
$z'\in{_b(\mathcal MU_c)}$.

\begin{theorem} \label{T:Action} Let $A$ be an algebra in
$\mathbf V$, and let $M$ be an $A$-module totally in
$\mathbf V$. Let
$Z$ be as just defined. We have
\begin{enumerate}
\item $(x)_bm=m$;
\item if $z\in{_aZ_b}$, then $z(x)_b=z$ and $(x)_az=z$;
\item the action $\pair zm\mapsto zm$ is bilinear;
\item the multiplication $\pair z{z'}\mapsto zz'$ is
bilinear;
\item if $u\in U$, and $m\in{_bM}$, then
$(u)_bm=u^M_{\langle b\rangle}(m)$;
\item if $u$, $v\in U$,
$(v)_{u(b)}(u)_b=(vu)_b$;
\item the action of $Z$ on $M$ is associative, i.e.,
$(z'z)m=z'(zm)$;
\item multiplication in $Z$ is associative;
\item $Z$ is a ringoid, with $1^Z_b=(x)_b$ and
$_a0^Z_b=(a)_b$; and
\item the $|A|$-tuple of abelian groups $M$, with the
defined action by
$Z$, form a left
$Z$-module.
\end{enumerate}
Now let $M'$ be another $A$-module totally in $\mathbf
V$, and $\phi:M\to M'$ a homomorphism. We have
\begin{enumerate}
\item[11.] If $z\in{_aZ_b}$ and
$m\in{_bM}$, then $z({_b\phi(m)})={_a\phi(zm)}$.
\end{enumerate}
\end{theorem}

\begin{proof} (1): This is clear from the definition of
the action of $Z$ on $M$.

(2) That $z(x)_b=z$ follows from the fact that
$\theta_{b,{(x)_b}}=1_{\mathcal MU_b}$. That $(x)_az=z$
follows from (1).

(3): Let $m\in{_bM}$ and $z\in{_aZ_b}$. If
$m=n_1m_1+n_2m_2$, then
$\theta_{b,m}=n_1\theta_{b,m_1}+n_2\theta_{b,m_2}$. Thus,
$zm={_a(\theta_{b,m})(z)}=n_1{_a(\theta_{b,m_1})(z)}
+n_2{_a(\theta_{b,m_2})(z)}=zm_1+zm_2$.
If $z=n_1z_1+n_1z_2$, then
$_a(\theta_{b,m})(z)=n_1{_a(\theta_{b,m})(z_1)}
+n_2{_a(\theta_{b,m})(z_2)}=z_1m+z_2m$.

(4): Follows from (3).

(5): We have
\begin{align*}
(u)_bm
&={_{u(b)}(\theta_{b,m})((u)_b)}\\
&={_{u(b)}(\theta_{b,m})(_{u(b)}\eta U_b)(u))}\\
&={_{u(b)}(\psi_{b,m})(u)}\\
&=u^{\mathcal UM}_{\langle b\rangle}(m)\\
&=u^M_{\langle b\rangle}(m),
\end{align*}
by the definition of $u^M_{\langle b\rangle}(m)$ from
\S\ref{S:Polynomials}.

(6): Using (5), we have
\begin{align*}
(v)_{u(b)}(u)_b
&= v^{Z_b}_{\langle u(b)\rangle}((u)_b)\\
&= v^{Z_b}_{\langle u(b)\rangle}({_{u(b)}(\eta U_b)(u)})\\
&={_{v(u(b))}(\eta U_b)(v^{U_b}_{\langle u(b)\rangle}(u))}
\\
&={_{v(u(b))}(\eta U_b)(v^{U,\iota_A}(u))}
\\
&={_{v(u(b))}(\eta U_b)(vu)}
\\
&=(vu)_b.
\end{align*}

(7): From theorem~\ref{T:LinearComb}, $z'$ and $z$ are
linear combinations respectively of elements of form
$(u')_{b'}$ and
$(u)_b$ where $b'=u(b)$. By (3), the statement to be
proved then reduces to the same statement for
$z=(u')_{u(b)}$ and
$z=(u)_b$. But this follows from (6), (5), and the
associativity of action for unary polynomials
(theorem~\ref{T:PolyAction}(1)).

(8): Follows from (7).

(9): Follows from (8), (4), and (2), and the fact that
by definition, $(a)_b$ is the zero
element of $_aZ_b$, because it is the image under
the pointed $A$-overalgebra homomorphism $_a(\eta U_b)$ of
$_a*^{U_b}$.

(10): Follows from (7), (3), and (1).

(11): Follows from the fact that
$\theta_{b,{_b\phi(m)}}=\phi\theta_{b,m}$.
\end{proof}

$Z$ is the \emph{enveloping
ringoid of
$A$ with respect to $\mathbf V$}. More formally, we
denote this ringoid by
$\Z[A,\mathbf V]$. However, we will continue to abbreviate
it by
$Z$ in what follows, when no confusion is likely.

\section{The Isomorphism of $\Z[A,\mathbf V]$-$\Mod$ and
$\Ab[A,\mathbf V]$}\label{S:Isomorphism}

From theorem~\ref{T:Action} (10) and (11), the
construction of the left
$Z$-module structure on $M$ is a functor
$\mathcal G:\Ab[A,\mathbf V]\to Z$-$\Mod$, where we
define
$\mathcal G\zeta=\zeta$ for $\zeta:M_1\to M_2$.

Given a left $Z$-module $M$, we place an $A$-module
structure $\mathcal HM$ on the $A$-tuple of abelian
groups
$_aM$ by defining
\[
\omega^{\mathcal HM}_{\mathbf
a,i}(m)=\left(\omega^U(a_1,\ldots,a_{i-1},x,
a_{i+1},\ldots,a_n)\right)_{a_i}(m),
\]
for each $n$-ary operation symbol $\omega$,
$n$-tuple $\mathbf a$, $i\in\{\,1,\ldots,n\,\}$, and
$m\in{_{a_i}M}$.

\begin{lemma} \label{T:LeadupLeadup}
If $\omega$ is an $n$-ary operation symbol,
$\mathbf a$ is an $n$-tuple of elements of $A$,
$1\leq i\leq n$, and $z\in{_{a_i}Z_b}$, then
$\omega^{Z_b}_{\mathbf a,i}(z)
=\left(\omega^U(a_1,\ldots,a_{i-1},
x,a_{i+1},\ldots,a_n)\right)_{a_i}z.$
\end{lemma}

\begin{proof} We have
\begin{align*}
\omega^{Z_b}_{\mathbf
a,i}(z)
&=\omega^{Z_b}_{\mathbf
a,i}({_{a_i}\theta_{a_i,z}(1_{a_i})})
\\
&={_{\omega(\mathbf
a)}\theta_{a_i,z}(\omega^{Z_{a_i}}_{\mathbf
a,i}(1_{a_i}))}
\\
&=\left(\omega^{Z_{a_i}}_{\mathbf a}({_{a_1}0},\ldots,
{_{a_{i-1}}0},x,{_{a_{i+1}}0},\ldots,{_{a_n}0})\right)z \\
&=\left({_{\omega(\mathbf a)}(\eta
U_{a_i})}\omega^{U_{a_i}} _{\mathbf
a}({_{a_1}*},\ldots,{_{a_{i-1}}*},x,{_{a_{i+1}}*},
\ldots,{_{a_n}*})\right)z \\
&=
\left(\omega^U(a_1,\ldots,a_{i-1},x,a_{i+1},\ldots,
a_n)\right)_{a_i}
z.
\end{align*}
\end{proof}

\begin{lemma} \label{T:InvLeadup} If
$t$ is an
$n$-ary term,
$\mathbf a$ an $n$-tuple of elements of
$A$, $i$ a number such that $1\leq i\leq
n$, and $m\in{_{a_i}M}$, then
\[t^{\mathcal HM}_{\mathbf
a,i}(m)=\left(t^U(a_1,\ldots,a_{i-1},
x,a_{i+1},\ldots,a_n)\right)_bm.
\]
\end{lemma}

\begin{proof}
Uses the fact that $\Clo_n\Omega$ is generated by the
$x_i$. For
$t=x_j$, we have
\[t^{\mathcal HM}_{\mathbf a,i}(m)
=\begin{cases}
{_{a_j}0}, & \text{if $i\neq j$, and } \\
m, & \text{if $i=j$.}
\end{cases}
\]
and
\begin{align*}
(t^U(a_1,\ldots,a_{i-1},x,&a_{i+1},\ldots,
a_n))_{a_i}m\\
&=\begin{cases}
(a_j)_{a_i}m, & \text{if $i\neq j$, and} \\
(x)_{a_i}m, & \text{if $i=j$.}
\end{cases}\\
&=\begin{cases}
{_{a_j}0},\quad & \text{if $i\neq j$, and} \\
1^Z_{a_i}m,\quad & \text{if $i=j$.}
\end{cases}
\end{align*}
Thus, the statement is true in that case.

Now let $\omega$ be an $k$-ary operation symbol, and let
the statement be true for a $k$-tuple $\mathbf s$ of
$n$-ary terms. That is, if for each $j$ we
write $u_j$ for $s_j^U(a_1,\ldots,a_{i-1},x,a_{i+1},
\ldots,a_n)$, we will be assuming $
(s_j)^{\mathcal HM}_{\mathbf a,i}(m)=(u_j)_{a_i}m$ for
all $m$. We
have
\begin{align*}
((\omega\mathbf s)^U(a_1,
&\ldots,a_{i-1},x,a_{i+1},\ldots,
a_n))_{a_i}m \\
&=\left(\omega^U(\ldots,u_j,
\ldots)\right)_{a_i}m \\
&=\left({_{\omega\mathbf s(\mathbf a)}(\eta U_{a_i})}
(\omega_{\mathbf s(\mathbf a)}^{U_{a_i}}(\ldots,
u_j,\ldots)
)\right)m \\
&=\left(\omega^{Z_{a_i}}_{\mathbf s(\mathbf a)}(\ldots,
(u_j)_{a_i},
\ldots,)\right)m \\
&=\left(\sum_j\omega^{Z_{a_i}}_{\mathbf s(\mathbf a),j}
(u_j)_{a_i}\right)m
\\
&=\left(\sum_j\left(\omega^U(s_1(\mathbf
a),\ldots,s_{j-1}(\mathbf a), x, s_{j+1}(\mathbf
a),\ldots,s_k(\mathbf a))\right)_{s_j(
\mathbf a)}
(u_j)_{a_i}\right)m
\\
&=\sum_j\left(\omega^U(s_1(\mathbf
a),\ldots,s_{j-1}(\mathbf a), x, s_{j+1}(\mathbf
a),\ldots,s_k(\mathbf a))\right)_{s_j(
\mathbf a)}
((u_j)_{a_i}m)
\\
&=\sum_j\omega^{\mathcal HM}_{\mathbf s(\mathbf
a),j}((u_j)_{a_i}m) \\
&=\omega^{\mathcal HM}_{\mathbf s(\mathbf
a)}\left(\ldots,(s_j)^{\mathcal HM}_{\mathbf
a,i}(m),\ldots\right)
\\
&=(\omega\mathbf
s)^{\mathcal HM}_{\mathbf s(\mathbf a),i}(m).
\end{align*}
\end{proof}

\begin{theorem} $\mathcal H$ is a functor from
$Z$-$\Mod$ to
$\Ab[A,\mathbf V]$, where for a homomorphism $\phi$,
$\mathcal H\phi=\phi$.
\end{theorem}

\begin{proof}
Let
$M$ and
$M'$ be left
$Z$-modules and
$\phi:M\to M'$ a homomorphism. To prove that $\mathcal
H$ is functorial, we must show that if
 $\omega$ is an $n$-ary operation symbol,
$\mathbf a\in A^n$, $1\leq i\leq n$, and
$m\in{_{a_i}M}$, then
\[\omega^{\mathcal HM}_{\mathbf a,i}({_{a_i}\phi(m)})=
{_{\omega(\mathbf a)}\phi(\omega^{\mathcal HM}_{\mathbf
a,i}(m))}.\]
Let us write $m'$ for $_{a_i}\phi(m)$. We have
\begin{align*}
\omega^{\mathcal HM}_{\mathbf a,i}({_{a_i}\phi(m)})
&=\left(\omega^U(a_1,\ldots,a_{i-1},x,a_{i+1},\ldots,
a_n)\right)_{a_i}
(m') \\
&=\left({_{\omega(\mathbf a)}(\eta U_{a_i})}
(\omega^U(a_1,\ldots,a_{i-1},x,a_{i+1},\ldots,a_n))
\right)(m') \\
&={_{\omega(\mathbf a)}(\theta_{a_i,m'})}
\left({_{\omega(\mathbf a)}(\eta
U_{a_i})(\omega^U(a_1,\ldots,
a_{i-1},x,a_{i+1},\ldots,a_n)})\right) \\
&=({_{\omega(\mathbf
a)}\phi})({_{\omega(\mathbf a)}
\theta_{a_i,m}})\left({_{\omega(\mathbf a)}(\eta U_i)}
(\omega^U(a_1,\ldots,a_{i-1},x,a_{i+1},\ldots,a_n))\right)
\\
&={_{\omega(\mathbf
a)}\phi\left(\left(\omega^U(a_1,\ldots,a_{i-1},
x,a_{i+1},\ldots,a_n)\right)_{a_i}(m)\right)} \\
&={_{\omega(\mathbf
a)}\phi(\omega^{\mathcal HM}_{\mathbf a,i}(m))}.
\end{align*}

Now let $M$ be a left
$Z$-module.
If $t=t'$ is an identity of $\mathbf V$, then by 
lemma~\ref{T:InvLeadup}, we have $t^{\mathcal HM}_{\mathbf
a,i}=(t')^{\mathcal HM}_{\mathbf a,i}$ for each $i$. 
Thus,
$\mathcal HM$ is totally in $\mathbf V$.
\end{proof}

The next theorem can be compared with the treatment
of
\cite[corollary 8.8]{R92}:

\begin{theorem}
The functors
$\mathcal G$ and $\mathcal H$ are inverses.
\end{theorem}

\begin{proof} Let $M$ be an $A$-module totally in
$\mathbf V$. We will show that for each $n$-ary operation
symbol $\omega$, $\omega^{\mathcal
H\mathcal GM}_{\mathbf a,i}=\omega^M_{\mathbf a,i}$,
for all $\mathbf a$ and $i$. Let
$m\in{_a(\mathcal H\mathcal GM)}={_a(\mathcal GM)}={_aM}$.
We have
\begin{align*}
\omega^{\mathcal H\mathcal GM}_{\mathbf a,i}(m)
&=\left(\omega^U(a_1,\ldots,a_{i-1},
x,a_{i+1},\ldots,a_n)\right)_{a_i}m \\
&=\left({_{\omega(\mathbf a)}(\eta U_{a_i})}
(\omega^{U_{a_i}}_{\mathbf a}(a_1,\ldots,a_{i-1},x,
a_{i+1},\ldots,a_n))\right)m \\
&=\left(\omega^{Z_{a_i}}_{\mathbf
a}((a_1)_{a_i},\ldots,(a_{i-1})_{a_i},
(x)_{a_i},(a_{i+1})_{a_i},\ldots,(a_n)_{a_i})\right)m \\
&=\theta^M_{a,m}(\omega^{Z_{a_i}}_{\mathbf a}(
(a_1)_{a_i},\ldots,(a_{i-1})_{a_i},(x)_{a_i},
(a_{i+1})_{a_i},\ldots,(a_n)_{a_i}))\\
&=\omega^M_{\mathbf a}({_{a_1}0_{a_i}}m,\ldots,
{_{a_{i-1}}0_{a_i}}m,1_{a_i}m,
{_{a_{i+1}}0_{a_i}}m,\ldots, {_{a_n}0_{a_i}}m) \\
&=\omega^M_{\mathbf a}({_{a_1}0},\ldots,
{_{a_{i-1}}0},m,{_{a_{i+1}}0},\ldots,{_{a_n}0}) \\
&=\omega^M_{\mathbf a,i}(m).
\end{align*}
Thus, $\mathcal H\mathcal G=1_{\Ab[A,\mathbf V]}$.

Now, let $M$ be a
left $Z$-module. We must show that for all $a$, $b\in A$,
for all $z\in{_aZ_b}$, and for all
$m\in{_b(\mathcal G\mathcal HM)}={_bM}$, we have the
action of
$z$ on
$m$ the same in $\mathcal G\mathcal HM$ as in $M$. It
suffices to show this for $z$ of the form $(u)_b$, since
such elements generate $_aZ_b$ as a group, and the
two actions are bilinear. Let $u=t^U(x,c_1,\ldots,c_n)$,
and denote
$\langle b, c_1,\ldots,c_n\rangle$ by $\langle b,\mathbf
c\rangle$. In $\mathcal G\mathcal HM$, we have
\begin{align*}
(u)_bm
&=\left(t^U(x,c_1,\ldots,c_n)\right)_bm \\
&=\left({_{u(b)}(\eta U_b)}(t^{U_b}_{\langle b,\mathbf
c\rangle}(x,c_1,\ldots, c_n))\right)m \\
&=\left(t^{Z_b}_{\langle b,\mathbf
c\rangle}(1^Z_b,{_{c_1}0_b},
\ldots,{_{c_n}0_b})\right)m \\
&=\theta^{\mathcal HM}_{b,m}\left(t^{Z_b}_{\langle b,
\mathbf c\rangle}(1^Z_b,{_{c_1}0_b},\ldots,{_{c_n}0_b})
\right)\\
&=t^{\mathcal HM}_{\langle b,\mathbf
c\rangle}(m,{_{c_1}0},\ldots, {_{c_n}0}) \\
&=t^{\mathcal HM}_{\langle b,\mathbf c\rangle,1}(m).
\end{align*}

However, this is $\left(t^U(x,c_1,\ldots,c_n)\right)_bm$
in
$M$, by lemma~\ref{T:InvLeadup}. Thus,
$\mathcal G\mathcal H=1_{Z\text{-}\Mod}$.
\end{proof}

\section{$\Z[A,\mathbf V]$ and $Z_M$}\label{S:Canonical}

If $M$ is an $A$-module totally in $\mathbf V$, then
the action on $Z$ on $M$ yields a group homomorphism
from $_aZ_b$ to $_a\End(M)_b=\Ab({_bM},{_aM})$ for
each $a$, $b\in A$. These mappings comprise a ringoid
homomorphism
$f_M:Z\to
\End(M)$, sending 
$(u)_b$ to $u^M_{\langle b\rangle}$.

\begin{theorem} We have
\begin{enumerate}
\item $f_M$ is a ringoid homomorphism;
\item $f_M$ has range in $Z_M$, i.e.,
$f_M(z)\in{_b(Z_M)_a}$ for all $z\in{_bZ_a}$; and
\item $f_M$, considered as a homomorphism from $Z$ to
$Z_M$, is
\emph{cofaithful} (i.e., one-one and onto on objects and
onto on all hom-sets).
\end{enumerate}
\end{theorem}

\begin{proof}
(1): It suffices to show that given $u$, $v\in U$, we have
$v^M_{\langle u(b)\rangle}u^M_{\langle
b\rangle}=(vu)^M_{\langle b\rangle}$. But this follows
from the fact that
$\Phi^M_{\Pol}$ is a clone homomorphism.

(2): $\Pol_1(A,\mathbf V)$ is generated by the constant
polynomials $a\in A$, and $x$. Thus we can prove that
$f_M((u)_b)\in Z_M$ for all $u$, by showing that the
generators satisfy this condition and that the subset of
$\Pol_1(A,\mathbf V)$ satisfying the condition is closed
under the basic operations. We have
$f_M((a)_b)=a^M_{\langle b\rangle}={_a0_b^{Z_M}}$ for all
$a$, and
$f_M((x)_b)=x^M_{\langle b\rangle}=1_{{_bM}}$. Finally,
given
$\omega\in\Omega_n$, and $\mathbf u\in\Pol_1(A,\mathbf
V)^n$, we have $f_M(\omega^{\Pol_1(A,\mathbf V)}(\mathbf
u))\in{_{\omega(\mathbf u(b))}(Z_M)_a}$. For,
$\omega^{\Pol_1(A,\mathbf V)}(\mathbf
u)=(\tilde\Phi(\omega))(\mathbf u)$ by
theorem~\ref{T:PolyDef}(4), so that we have
\begin{align*}
f_M((\omega^{\Pol_1(A,\mathbf V)}(\mathbf u))_b)(m)
&=f_M((\tilde\Phi(\omega)(\mathbf u))_b(m)\\
&=((\tilde\Phi(\omega)(\mathbf u))^M_{\langle
b\rangle}(m)\\
&=((\tilde\Phi(\omega))^M\circ_{\mathbf
u^A}\mathbf u^M)_{\langle b\rangle}(m)\\
&=(\tilde\Phi(\omega))^M_{\mathbf u(b)}
\left((u_1)^M_{\langle
b\rangle}(m),\ldots,(u_n)^M_{\langle b\rangle}(m)\right)\\
&=\sum_i\left(\omega^M_{\mathbf
u(b),i}\right)\left((u_i)^M_{\langle b\rangle}(m)\right).
\end{align*}

(3): Obviously, $f_M$ is one-one and onto on objects. Let
$\omega\in\Omega_n$, and $\mathbf a\in A^n$, and consider
$\omega^M_{\mathbf a,i}\in{_{\omega(\mathbf
a)}(Z_M)_{a_i}}$. Let $\mathbf w=\langle
a_1^A,\ldots,a_{i-1}^A,x^A,a_{i-1}^A,\ldots,a_n^A\rangle
\in(\Clo_1|A|)^n$.
We have
\begin{align*}
f_M\left(\left(\omega(a_1,\ldots,a_{i-1},x,a_{i+1},
\ldots,a_n)\right)
_{a_i}\right)&(m)\\
=&\left(\omega(a_1,\ldots,a_{i-1},x,a_{i+1},
\ldots,a_n)\right)^M_{
\langle a_i\rangle}(m)\\
=&(\omega^M\circ_{\mathbf w}\langle a_1^M,\ldots,
a_{i-1}^M,x^M,a_{i+1}^M,\ldots,a_n^M\rangle)_{\langle
a_i\rangle}(m)\\ =&\omega^M_{\mathbf a}({_{a_1}0},\ldots,
{_{a_{i-1}}0},m,{_{a_{i+1}}0},\ldots,{_{a_n}0})\\
=&\omega^M_{\mathbf a,i}(m).
\end{align*}
The range of $f_M$ therefore contains all the generators
of $Z_M$, and thus is onto on all
hom-sets.
\end{proof}

\section{$J$ and $R$}\label{S:JandR}

For every $a$, $b\in A$, we have $_aZ_b={_a(\hat
\mathcal MU_b/\mathcal KU_b)}$. Define $_a\hat
Z_b={_a(\hat\mathcal MU_b)}$. It is easy to define a
ringoid structure on
$\hat Z$ by extending the composition of unary
polynomials $\Z$-bilinearly.  The natural maps from
$_a\hat Z_b$ to $_aZ_b$ are abelian group
homomorphisms and by theorem~\ref{T:Action}(6), constitute
a ringoid homomorphism. Thus, the kernel groups
$_aJ_b={_a(\mathcal KU_b)}$ constitute a ringoid ideal $J$
of
$\hat Z$.

We want to give another
description of the ideal
$J$, by giving generators for each abelian group $_aJ_b$.
For each
$n\geq 0$, each
$n$-ary polynomial
$\Pi$, and each
$b\in A$, we define
\[
R_{\Pi,b}=\Pi^{U,\myiota}(x,\ldots,x)
-\sum_{j=1}^n
\Pi^{U,\myiota}
(b,\ldots,b,x,b,\ldots,b),
\]
where in the sum over $j$, the $j^{\text{th}}$ argument
of $\Pi^{U,\myiota}$ is
$x$. For each $a$, $b\in A$, let $_aR_b$ be the subgroup
of
$_a(\hat\mathcal MU_b)$ generated by the
elements
$R_{\Pi,b}$, where
$\Pi(b,\ldots,b)=a$.

\begin{theorem}
For all $a$, $b\in A$, $_aR_b={_aJ_b}$.
\end{theorem}

\begin{proof}
Let $\Pi$ be an $n$-ary polynomial
such that $\Pi(b,\ldots,b)=a$. We will
show that
$R_{\Pi,b}\in{_aJ_b}$.

First, let $\mathbf c$ be an $m$-tuple of elements of
$A$, and $t$ an $(n+m)$-ary term, such that
$\Pi=t^{\Pol_nA}(x_1,\ldots,x_n,c_1,\ldots,c_m)$. We have
\[t^U(x,\ldots,x,c_1,\ldots,c_m)=t^{U_b}_{\langle
b,\ldots,b,c_1,\ldots,c_m\rangle}(x,\ldots,x,c_1,\ldots,
c_m).\] Also, by
\eqref{E:Termfact}, we have for each
$i\in\{\,1,\ldots,n\,\}$,
\[
t^{\hat\mathcal MU_b}_{\langle
b,\ldots,b,c_1,\ldots,c_m\rangle,i}(x)- t^{U_b}_{
\langle
b,\ldots,b,c_1,\ldots,c_m\rangle}(b,\ldots,b,x,b,\ldots,
b,c_1,\ldots,c_m)\in{_aJ_b},
\]
where in the second term, $x$ appears as the
$i^{\text{th}}$ argument.
Finally, we have
\[
t^{U_b}_{\langle b,\ldots,b,c_1,\ldots,c_m\rangle}
(x,\ldots,x,c_1,\ldots,c_m)
-\sum_{i=1}^n
t^{\hat\mathcal MU_b}
_{\langle b,\ldots,b,c_1,\ldots,c_m\rangle,i}
(x)\in{_aJ_b},
\]
because the $a$-component of the
natural map from $\hat\mathcal MU_b$ to $\mathcal
MU_b$ sends this element to
\[
t^{\mathcal
MU_b}_{\langle b,\ldots,b,c_1,\ldots,c_m\rangle}
(x/K,\ldots,x/K,{_{c_1}0},\ldots,{_{c_m}0})
-\sum_{i=1}^n
t^{\mathcal MU_b}_{\langle
b,\ldots,b,c_1,\ldots,c_m\rangle,i}(x/K),\] where
$K$ stands for $_b(\mathcal KU_b)$, and this
is
$_a0$.

Thus, we have
\begin{align*}
R_{\Pi,b} &= \Pi^{
U,\myiota}(x,\ldots,x)-\sum_{i=1}^n
\Pi^{U,\myiota}(b,\ldots,b,x,b,\ldots,b) \\
&=t^U(x,\ldots,x,c_1,\ldots,c_m)-\sum_{i=1}^n
t^U(b,\ldots,b,x,b,\ldots,b,c_1,\ldots,c_m) \\
&=t^{U_b}_{\langle b,\ldots,b,c_1,\ldots,c_m\rangle}
(x,\ldots,x,c_1,\ldots,c_m)-\sum_{i=1}^n
(t^{\hat\mathcal MU_b}
_{\langle b,\ldots,b,c_1,\ldots,c_m\rangle,i}(x)+k_i)\\
&\in{_aJ_b},
\end{align*}
where $k_i\in{_aJ_b}$ for each $i$.

Since the abelian group generators of $_aR_b$ belong to
the abelian group $_aJ_b$, we have
$_aR_b\subseteq{_aJ_b}$.

To show $_aJ_b\subseteq{_aR_b}$, we will show that for
each $b\in A$, the elements of $_a\mathcal SU_b$ (where
$\mathcal SU_b$ is the $A$-set defined in
section~\ref{S:Modulization} that generates
$\mathcal KU_b$) are contained in
$_aR_b$, and also that for each
$b$, the $A$-tuple of abelian groups $_aR_b$ is closed
under action by elements of $Z_{\hat\mathcal MU_b}$, i.e.,
that the
$_aR_b$ form an $A$-submodule of $\hat Z_b$.

We start with the fact that
$_a*^{U_b}\in{_aR_b}$ for each $a\in A$. For,
$_a*^{U_b}=a\in\Pol_0A$. Next, let
$\omega$ be an $n$-ary operation symbol, $\mathbf
a$ an
$n$-tuple of elements of $A$, and $\mathbf u\in{_{\mathbf
a}U_b}$. We need
to show that $\omega^{U_b}_{\mathbf a}(\mathbf
u)-\sum_{i=1}^n\omega^{\hat\mathcal MU_b}_{\mathbf
a,i}(u_i)$ belongs to $_{\omega(\mathbf a)}R_b$. Now, each
$u_i$ is an element of $_{a_i}U_b$, which is a subset
of $U$, an algebra
generated by $x$ and elements of $A$. Thus, there are an
$m$-tuple $\mathbf d$ of elements of $A$, and $(m+1)$-ary
terms $s_i$, so that for each $i$ we have
\[u_i=s_i^U(d_1,\ldots,d_m,x).
\]
For each
$i$, define the $n$-ary polynomial
\[\Pi_i=s_i(d_1,\ldots,d_m,x_i),\]
and let $\Pi'=\omega^{\Pol_nA}(\Pi_1,\ldots,\Pi_n)$.
Then
we have
\[\omega^{U_b}_{\mathbf a}(\mathbf u)
-\sum_{i=1}^n\omega^{\hat\mathcal MU_b}_{\mathbf
a,j}(u_i) =R_{\Pi',b}
\in{_{\omega(\mathbf a)}R_b},
\]
proving that $\mathcal SU_b\subseteq R_b$.

Now we must show that homomorphisms of the form
$\omega^{\hat\mathcal MU_b}_{\mathbf c,i}$, where
$\omega$ is $n$-ary and $c_i=a$, send generators of
$_aR_b$ to elements of $_{\omega(\mathbf
c)}R_b$.
 Consider the generator
$R_{\Pi,b}$, where $\Pi$ is
$m$-ary and $\Pi(b,\ldots,b)=a=c_i$. We have
\begin{align*}
\omega^{\hat\mathcal MU_b}_{\mathbf c,i}(R_{\Pi,b})
&=\omega^U(c_1,\ldots,c_{i-1},\Pi^{
U,\myiota}(x,\ldots,x),c_{i+1},
\ldots,
c_n)\\
&\qquad-\sum_{j=1}^m\omega^U(c_1,\ldots,c_{i-1},
\Pi^{U,\myiota}(b,\ldots,b,x,b,\ldots,b),
c_{i+1},\ldots,c_n)\\
&=R_{\bar\Pi,b},
\end{align*}
where in the summation, $x$ appears in the
$j^{\text{th}}$ position, and $\bar\Pi$ is the
$m$-ary polynomial
\[\omega^{\Pol_m(A,\mathbf V)}(c_1,\ldots,c_{i-1},\Pi,
c_{i+1},\ldots,c_n).
\]
Thus, $_aJ_b\subseteq{_aR_b}$.
\end{proof}

\section{Previous Constructions of Enveloping Ringoids}
\label{S:Previous}

The second construction of the enveloping ringoid in
\cite{R92}, in the notation of this paper, defined
$\Z[A,\mathbf V]$ as
$\hat Z/R$, with only a slight difference in the
description of the ideal
$R$. Thus, the enveloping ringoid of this paper is the
same as the one defined in
\cite{R92}.

\section{Congruence-Modular Varieties} \label{S:CM}

For congruence-modular $\mathbf V$, then results
applying to that case
\cite{R92,R96} yield some simplifications
which shed light on the modulization functor and
are helpful for finding the structure of the enveloping
ringoid:

\begin{theorem} [\cite{R92}, Theorem~C.10.4] Let $P$ be a
pointed
$A$-overalgebra totally in $\mathbf V$, a
congruence-modular variety. Let $\kappa_P$ denote the
kernel congruence of the projection homomorphism
$\pi_P:A\semitimes P\to A$, and let
$\kappa$ denote $[\kappa_P,\kappa_P]$. Then
$\nat\kappa_P$ factors as $\pi\circ\nat\kappa$ where
$\pi:(A\semitimes P)/\kappa\to A$ is an onto
homomorphism. Let $\iota=(\nat\kappa)\iota_P$. Then
$\pi\iota=1_A$. Let
$M$ denote
$\lsil (A\semitimes P)/\kappa,\pi,\iota\rsil$. Then
$\pair M\phi$ is a universal arrow from
$P$ to $\mathcal U$, where $\phi:P\to\mathcal UM$ is
defined by $_a\phi:p\mapsto\pair ap/\kappa$.
\end{theorem}

\begin{corollary} Let $P$ be a pointed $A$-overalgebra
totally in $\mathbf V$, a congruence-modular variety.
Then $\eta P$ is onto. I.e., for each $a\in A$, every
element of $_a\mathcal MP$ is of the form $_a(\eta P)(p)$
for some $p\in{_aP}$.
\end{corollary}

Thus, when $\mathbf V$ is congruence-modular, every
element of
$_a\Z[A,\mathbf V]_b$ is of the form $(u)_b$ for some
$u\in\Pol_1(A,\mathbf V)$.

Another consequence of this corollary is the fact that if
$P$ is a pointed $A$-overalgebra totally in $\mathbf V$
and $\mathbf V$ is congruence-modular, then $\mathcal MP$
is totally in $\mathbf V$. Of course, we proved this more
generally for all $\mathbf V$ in \S\ref{S:Modulization},
but it required more effort to prove without the
assumption that $\mathbf V$ is congruence-modular.

\begin{theorem} Let $A$ be an algebra in $\mathbf V$, a
congruence-modular variety, and let $d$ be a ternary
difference term for $\mathbf V$. Let $M$ be an
$A$-module totally in $\mathbf V$. Then for each $a\in
A$, and $m$, $m'$, $m''\in{_aM}$ we have
\[
m-m'+m''=d^M(m,m',m'').
\]
\end{theorem}

\begin{corollary}
Let $A$ be an algebra in $\mathbf V$, a
congruence-modular variety, and let $d$ be a ternary
difference term for $\mathbf V$.  If
$b\in A$ and $u$, $u'$, $u''\in\Pol_1(A,\mathbf V)$ are
such that $u(b)=u'(b)=u''(b)=a$, then
$(d^U(u,u',u''))_b=(u)_b-(u')_b+(u'')_b$.
\end{corollary}

\begin{proof} We have 
\begin{align*}
(d^U(u,u',u''))_b
&=
(d^{U_b}_{\langle
a,a,a\rangle}(u,u',u''))_b \\
&=
{_a(\eta U_b)(d^{U_b}_{\langle a,a,a\rangle}(u,u',u''))}
\\
&=
d^{Z_b}_{\langle a,a,a\rangle}((u)_b,(u')_b,(u'')_b) \\
&= (u)_b-(u')_b+(u'')_b.
\end{align*}
\end{proof}

\subsection{Underlying groups}

In many cases of interest, a variety $\mathbf V$ has an
underlying group functor. In such a case, given
$A\in\mathbf V$ and $a$, $b\in A$, let $u=xa^{-1}b$ and
$v=xb^{-1}a$. Then we have $uv=x$ and $vu=x$. Thus,
$(u)_{v(b)}(v)_b={_a1^Z_b}$, and
$(v)_{u(a)}(u)_a={_b1^Z_a}$. It follows that all the
endomorphism rings $_aZ_a$ are isomorphic, and we have
situations like those listed in table~1. See \cite{R96}
for a more extensive discussion.

\end{document}